\documentclass[review,onefignum,onetabnum]{siamart190516}

\newcommand{\review}[1]{#1}
\newcommand{\rreview}[1]{#1}
\newcommand{\rrreview}[1]{{\color{black}{#1}}}

\usepackage{lipsum}
\usepackage{amsfonts, amssymb, dsfont}
\usepackage{graphicx}
\usepackage{epstopdf}
\usepackage{algorithmic}
\ifpdf
  \DeclareGraphicsExtensions{.eps,.pdf,.png,.jpg}
\else
  \DeclareGraphicsExtensions{.eps}
\fi

\newsiamremark{remark}{Remark}
\newsiamremark{hypothesis}{Hypothesis}
\crefname{hypothesis}{Hypothesis}{Hypotheses}
\newsiamthm{claim}{Claim}

\headers{State estimation from elastography images}{F. Galarce \textsl{et. al.}}

\title{
Displacement and pressure reconstruction from magnetic resonance elastography images: application to an in silico brain model
}
\author{
Felipe Galarce\footnote{S\lowercase{chool of} C\lowercase{ivil} E\lowercase{ngineering}, P\lowercase{ontificia} U\lowercase{niversidad} C\lowercase{at\'olica de} V\lowercase{alpara\'iso}, C\lowercase{hile}}\hspace{.05cm} $^\dagger$
\and Karsten Tabelow\thanks{Weierstrass-Institut f\"ur Angewandte Analysis und Stochastik. Leibniz-Institut im Forschungsverbund Berlin e.~V., Berlin, Germany}
\and Jörg Polzehl$^\dagger$
\and Christos Panagiotis Papanikas\thanks{Department of Mechanical \& Manufacturing Engineering, University of Cyprus, Nicosia, Cyprus.} 
\and Vasileios Vavourakis$^\ddagger$\thanks{Department of Medical Physics \& Biomedical Engineering, University College London, London, UK.} 
\and Ledia Lilaj\thanks{Department of Radiology, Charité–Universitätsmedizin Berlin, Campus Charité Mitte, Charitéplatz. 1, Berlin, Germany.}
\and Ingolf Sack$^\P$
\and Alfonso Caiazzo$^\dagger$
}
\usepackage{amsopn}
\usepackage{mathtools}
\usepackage{subfigure}
\usepackage[left=0.9in,right=0.9in,top=1in,bottom=1in,includefoot,includehead,headheight=13.6pt]{geometry}
\usepackage{nicefrac}
\definecolor{BlueMAP}{rgb}{0.27,0.45,0.77}

\newcommand{\cA}{\ensuremath{\mathcal{A}}}

\newcommand{\cL}{\ensuremath{\mathcal{L}}}
\newcommand{\cM}{\ensuremath{\mathcal{M}}}
\newcommand{\cN}{\ensuremath{\mathcal{N}}}

\newcommand{\cP}{\ensuremath{\mathcal{P}}}

\newcommand{\cW}{\ensuremath{\mathcal{W}}}

\newcommand{\bI}{\ensuremath{\mathbb{I}}}

\newcommand{\bP}{\ensuremath{\mathbb{P}}}

\newcommand{\bR}{\ensuremath{\mathbb{R}}}

\newcommand{\bg}{\ensuremath{\textbf{g}}}
\newcommand{\bn}{\ensuremath{\textbf{n}}}
\newcommand{\bv}{\ensuremath{\textbf{v}}}

\newcommand{\bz}{\ensuremath{\textbf{z}}}
\newcommand{\bu}{\ensuremath{\textbf{u}}}

\def\[{\left[}
\def\]{\right]}
\def\<{\langle}
\def\>{\rangle}
\def\({\left(}
\def\){\right)}
\def\[{\left [}
\def\]{\right]}
\def\({\left(}
\def\){\right)}

\newcommand{\norm}[1]{\Vert #1 \Vert}

\DeclareMathOperator*{\arginf}{arg\,inf}

\newcommand{\dx}{\ensuremath{\mathrm dx}}

\newcommand{\dt}{\ensuremath{\mathrm dt}}

\usepackage{soul}

\definecolor{revision}{HTML}{16A085}

\newtheorem{problem}[theorem]{Problem}

\newcommand{\manit}{\ensuremath{\cM^{\text{training}}}}

\begin{document}

\maketitle

\begin{abstract}
Magnetic resonance elastography is a motion-sensitive image modality that allows
to measure in vivo tissue displacement fields in response to mechanical excitations.
This paper investigates a data assimilation approach for reconstructing tissue displacement and pressure
fields in an in silico brain model from partial elastography data.
The data assimilation is based on a parametrized-background data weak methodology, in which the state of the physical system -- tissue displacements and pressure fields -- is reconstructed from the available data assuming an underlying poroelastic biomechanics model.
For this purpose, a physics-informed manifold is built by sampling the space of parameters describing the tissue model close to their physiological ranges to simulate the corresponding poroelastic problem, and computing a reduced basis via Proper Orthogonal Decomposition.
Displacements and pressure reconstruction is sought in a reduced space after solving a minimization problem that encompasses both the structure of the reduced-order model and the available measurements.
The proposed pipeline is validated using synthetic data obtained after simulating the poroelastic mechanics of a physiological brain.
The numerical experiments demonstrate that the framework can exhibit accurate joint reconstructions of both displacement and pressure fields.
The methodology can be formulated for an arbitrary resolution of available displacement data from pertinent images.
It can also inherently handle uncertainty on the physical parameters of the mechanical model by enlarging the physics-informed manifold accordingly.
Moreover, the framework can be used to characterize, in silico, biomarkers for pathological conditions by appropriately training the reduced-order model.
A first application for the non-invasive estimation of ventricular pressure as an indicator of abnormal intracranial pressure is shown in this contribution.
\end{abstract}



\maketitle

\section{Introduction}

Medical imaging, combined with in silico models and computer simulations, has an enormous potential to provide clinically valuable insights into the mechanics of complex biological tissues and to support noninvasive diagnostics.
In this context, data assimilation methods are being developed to complement the available data
with physical models, advanced numerical simulations, and mathematical methods.
This paper focuses on the assimilation of internal displacement data, as those acquired via magnetic resonance elastography (MRE) imaging, into a human brain in silico model.

MRE is a \textit{tissue} imaging modality designed to measure mechanical properties of biological tissues. 
It combines phase-contrast MRI with the propagation of harmonic mechanical waves (10--100 Hz) induced by actuators placed externally on the patient body surface \cite{muthupillai_1995}.
The mechanical response of the tissue is recorded as a three-dimensional \textit{internal} displacement field.
Combined with physical tissue models and different inversion methods, these displacement data allow obtaining quantitative information on tissue mechanical properties, e.g., in terms of mechanical parameters (see \cite{fovargue_2018_review, sack-2008, sack-bioqic-18}).

The clinical potential of elastography has constantly been increasing in the last decades \cite{manduca_etal_2021_review}.
It has been proven to be an effective approach for the quantitative estimation of biomarkers (such as elastic parameters, tissue fluidity, viscoelasticity) related to different tissue pathologies.
Relevant clinical applications of MRE include diagnosing and staging diseases that directly influence tissue stiffness, such as cancer and fibrosis \cite{BERTALAN2020395,Reiter_etal_2021,Streitberger_etal_2020}. 
Focusing on the brain, elastography has been used for the characterization of cancer tissue \cite{bunecicius_etal_2020} and for the early-stage diagnosis of neurological diseases characterized by alteration of microstructure properties of brain tissue (see, e.g.,  \cite{hiscox_etal_brain_mre_2016}).

This work is motivated by the applicability of elastography in monitoring and quantifying the increase of intracranial pressure (ICP).
Intracranial hypertension might be responsible for different neurological diseases, such as ischemia, tumors, and hydrocephalus \cite{ren_etal_2020}, as well as neurological disturbances and cerebral damage.
Recent research focused on the correlation between hydrocephalus and changes in elastic behavior (see \cite{fattahi_etal_2016} or the recent review \cite{aunan-diop_etal_2022}). 
In \cite{kreft_etal_2022} a novel method based on ultrasound time-harmonic elastography has been described, in which the quantification of shear wave speed is used as indirect biomarker for elevated ICP.

\review{This paper proposes a novel data assimilation pipeline used to reconstruct an approximation of the displacement and pressure solutions over the whole brain from partial displacement observations. A primary objective of this approach is to enable, }
for the first time, a computational method for non-invasive quantification of elevated ICP from partial displacement data as those acquired in a typical MRE examination.
Currently, precise measurement of ventricular pressure and ICP are based on invasive procedures such as catheterization or perforation and, hence, not suited for early-stage risk quantification and constant monitoring.
Non-invasive alternatives for estimating relevant pressure gradients would consequently drastically enhance the possibilities for diagnosis and patient monitoring.
In order to characterize the fluid pressure field within the brain tissue, it is necessary to consider a biphasic tissue model describing both the solid and the fluid (i.e., the ICP) mechanics.

Recent experimental results suggested that accounting for the fluid phase, MRE can be utilized to infer the presence of pathological pressure conditions by measuring the effect of the increase in interstitial pressure on tissue mechanics \cite{hirsch-etal-2014-liver,Hetzer-18,solamen_etal_2021}.
More recent studies address different applications of MRE with a biphasic model of the tissue, e.g., a poroelastic medium, in which the solid and interstitial fluid interactions have been investigated \cite{mcgarry_etal_2019, lilaj_etal_2021a,leiderman_etal_2006, mcgarry_etal_2015, mcgarry_etal_2019, pattison_etal_2014}.

Computational modeling of brain biomechanics is a challenging problem due to the structural complexity of the brain and the lack of suitable experimental data for model parametrization.
Biomechanical models of the human brain based on in vitro mechanical experiments have been recently studied in \cite{budday_etal_2017}, while a detailed discussion of suitable constitutive models for brain tissue based on in silico experiments has been presented in \cite{de_rooji_kuhl_2016}.
Focusing on hydrocephalus, the mechanics of ventricle growth has been computationally investigated in \cite{dutta-roy_etal_2008} using a realistic three-dimensional single-phase brain model based on hyperelastic constitutive law.
Biphasic (poro-viscoelastic) models have been investigated in \cite{comellas_etal_2020}, performing load experiments in silico to computationally investigate the response of the tissue.
The choice of the brain mechanical model shall be dictated by the regime of dynamics of interest and the targeted application \cite{budday_etal_2020a}.
While hyperelastic laws are suited for capturing long-term dynamics, this paper focuses on a linear poroelastic description. %
\review{The poroelastic model 
explicitly encompasses the fluid pressure  
as part of the solution state. At the same time,
the model is assumed to be able to describe the mechanical response to 
the vibration induced during MRE examination, see for example the  recent experimental and numerical works on MRE 
\cite{hirschbraunsack-17,perrinez_etal_2010,Tan-etal-2017}}.
Despite the linearity of the model, this problem is \rreview{highly} demanding for two major reasons: (i) the high dimension of the unknowns (the pressure field, defined at each point of a suitable discretization of the brain) and (ii) the limited availability and resolution of images.
In fact, MRE acquisition is practically constrained by the length of examination time, and displacement data are typically available only on a sub-region of the tissue of interest. 
In addition, MRE measurement are only available for the displacement field.

A finite element method for reconstructing poroelastic parameters from MRE 
(MR poroelastography) was first presented in \cite{perrinez_etal_2009, perrinez_etal_2010}, by considering a synthetic phantom. A numerical framework for the estimation of poroelastic interstitial pressure from MRE data has been recently proposed in \cite{Tan-etal-2017}.
The method is based on solving Biot problems locally in different patches located in the region where MRE data are available. 
To overcome the lack of pressure boundary conditions in the measurement regions, displacement data are used to derive suitable expressions for the pressure on the boundary.
The results of Tan \textit{et al.} \cite{Tan-etal-2017} are limited to a numerical phantom (cubic tissue sample with a single inclusion).
However, their approach still faces several challenges related to  robustness (noisy data) and due to the high number of unknown parameters, including the internal pressure boundary conditions between the patches.
A different approach for the estimation of the pressure field in the context of MRE has been recently proposed by Fovargue and colleagues \cite{fovargue-etal-2018, fovargue-etal-2020}, combining a stiffness reconstruction approach (also for large strains) with an analytic model of an inflating sphere to relate it to the pressure field.

The overarching goal of this paper is to assimilate MRE data of poroelastic biological tissues into in silico modeling from the perspective that in clinical applications, only \rreview{partially available displacement} measurements -- few slices and with limited resolution -- are available.
Thus, focusing here on brain elastography, \textit{reconstruction} of a poroelastic solution (both displacement and pressure fields) over the whole computational domain, the brain organ, is challenges. 
Full organ reconstruction is vital to infer pressure-dependent biomarkers if the regions of interest are not fully contained in the imaged subdomain.
This is the case, for instance, of the ventricular CSF pressure, which can be considered as an indicator of hydrocephalus.

\begin{figure}[!ht]
\centering
\includegraphics[width=0.8\textwidth]{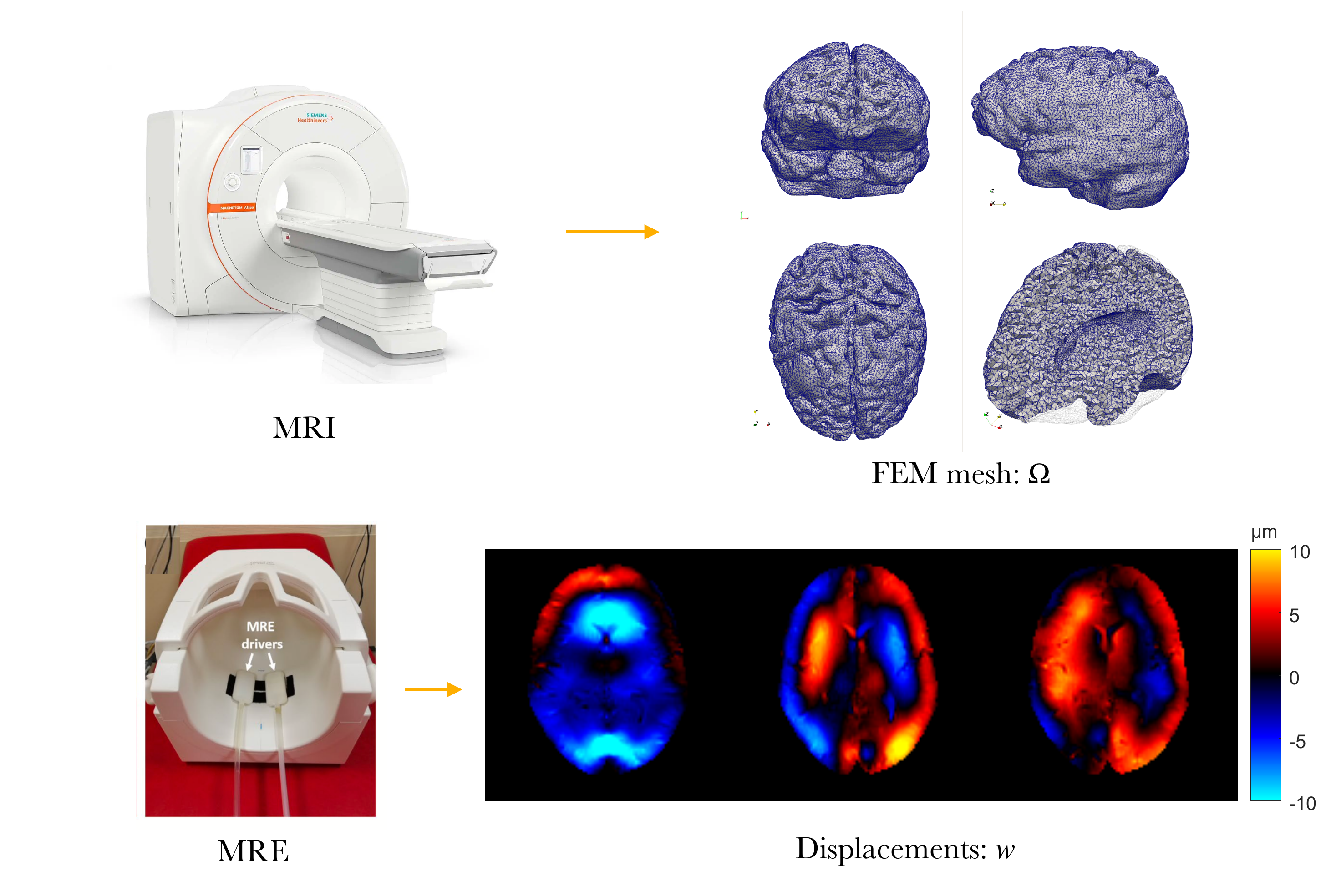}
\caption{The purpose of the considered data assimilation algorithm is to map partially available MRE measurements acquired on a few slices onto a displacement field defined on a three-dimensional brain model. As an example, a typical wave image on a slice is shown in the bottom sketch. The mapping requires anatomical data to construct the brain model, and it is based on assuming an underlying tissue model (encoded in a PDE), which is numerically solved via a finite element method to construct a suitable space for the reconstruction (top sketch).}
\label{fig:data_assimilation_sketch}
\end{figure}

The considered reconstruction algorithm addresses the state estimation problem based on the parametrized-background data-weak (PBDW) method \cite{maday-pbdw-2014}.
Originally introduced for wave equations, the PBDW has been recently extended in the context of ultrasound imaging and PC-MRI images of blood flow \cite{GGLM2021,galarceThesis,GLM2022}.
This approach is designed to \review{address the reconstruction of 
a physical solution over the whole domain when only partial displacements information
is available.
To this purpose, the reconstruction algorithm is based on solving a background parametrized PDE
on a patient-specific finite element model.}
The underlying physics is then exploited to extend the available measurement to a function defined in the whole three-dimensional brain model (Figure \ref{fig:data_assimilation_sketch}).
\review{In the considered case, the computational model has been 
created from three-dimensional anatomical MRI data acquired from the same subject.
Recent research has addressed the possibility of mapping geometrical information 
from existing models \cite{guibert-etal-2014,GLM2022} to further reduce the amount of data required for the reconstruction.}

To solve the problem, one generates a PDE-informed manifold of suitable solutions of the underlying equation -- a training set -- sampling the parameter space. 
\review{The reconstruction of the full state is computed by solving a constrained optimization problem
on a low-dimensional subspace, constructed by combining a reduced-order model 
(proper orthogonal decomposition) of the PDE-informed manifold and the linear functionals describing the available observations.
In this respect, the PBDW differs from standard approaches based, e.g., on Tikhonov regularizations or 
stabilization of the resulting finite element formulation \cite{burman-2020-1}.}
The method can be interpreted as a predictor-corrector scheme, in which a solution belonging to the low-dimensional subspace of the training set is improved by adding a correction orthogonal to the subspace, that accounts for the model bias \cite{GGLM2021,maday-pbdw-2014}.
By construction, the method avoids the issue of setting internal pressure boundaries and can naturally take into account a small variability in the physical and mechanical parameters.

The main contribution of this work is twofold:
First, a detailed validation of PBDW in the context of brain elastography is presented.
The validation is based on a pipeline with multidisciplinary components, entailing multi-modal image acquisition, data processing and image segmentation, in silico model generation, numerical solution, and data assimilation.
In this paper, the assimilated displacement data are limited to synthetic images.
However, the computational setup considers a physiological brain model, and the synthetic data have been defined in such a way to mimic a realistic experimental setup used in recent clinical research works \cite{lilaj_etal_2021a}. 
Second, the framework is used for the characterization of pressure-dependent biomarkers, such as ventricular pressure, to characterize pathological ICP gradients.
In particular, this characterization is done using only elastography data on a subset of the domain, not necessarily including the ventricles.
For the present application demonstration, two populations of synthetic patients are created, classified into a healthy and an unhealthy group.
The pathological states are characterized by a pressure increase in the ventricle by 10\% with respect to the healthy ones. 
For both cases, forward numerical simulations are performed, sampling displacement data on a few slices at the level of the eyes.
By properly constructing the training set for the reconstruction, the data assimilation algorithm is able to correctly classify all cases in the corresponding group.
The application of the proposed framework to in vivo displacement data is currently under investigation, combined with novel pre-processing techniques to estimate tissue properties from MRI data \cite{lilaj_etal_2021a, lilaj_etal_2021b}.

The rest of the paper is organized as follows.
Section \ref{sec:pbdw} introduces the PBDW data assimilation framework from a general perspective, while Section \ref{sec:forward} describes the considered poroelastic model and the numerical method used for the solution.
The application of the PBDW to the particular context is detailed in Section \ref{sec:mre-assimilation}.
The validation in the case of synthetic MRE measurements and the application to hydrocephalus classification are shown in Section \ref{sec:results}, while the conclusions are drawn in Section \ref{sec:conclusions}.

\section{A framework for data assimilation}\label{sec:pbdw}

The proposed numerical framework is based on the Parametri\-zed
Background Data-Weak (PBDW) method for state estimation, originally proposed in \cite{MPPY2015}, and recently applied also in the context of 
cardiovascular imaging \cite{GGLM2021, GLM2021}.

Let us consider a physical domain $\Omega \subset \bR^d$ ($d=2,3$), and let us assume that tissue biomechanics is described by a parameter-dependent partial differential equation $\cP$ of the form

\begin{equation}
\cP(u,\theta) = 0
\,,
\label{eq:model-pde}
\end{equation}

\noindent
which will be referred thereafter as the \textit{forward model}.
In \eqref{eq:model-pde}, $u \in V$ denotes the \textit{state}, and $\theta \in \bR^P$ stands for an array of \textit{parameters}. 

The rest of this section describes
the pipeline for the data assimilation of a general forward model and general observations, introducing the main notations and the different computational steps.
In Section \ref{sec:forward}, it will be shown how to apply the framework to the particular case of a forward model based on the Biot equations -- the state $u$ consists of displacement and pressure fields -- and to the case of observations defined by MRE data.

\subsection{PBDW preliminaries} 

The space $V$ (also called the \textit{ambient space}) is a functional space suitable for obtaining a well-posed forward problem \eqref{eq:model-pde}. It will be assumed to be a Hilbert space, endowed with an inner product, which will be denoted by $\<\cdot, \cdot\>$, and with the induced norm $\norm{\cdot} = \sqrt{\<\cdot , \cdot\>}$.
Depending on the PDE under consideration, typical choices for $V$ are subspaces of $L^2(\Omega)$ or $H^1(\Omega)$ (or product spaces of those), with the corresponding standard inner products.

Moreover, let us consider an appropriate discretization of the forward model \eqref{eq:model-pde} based on a finite element method.
As such, let us introduce a computational mesh $\mathcal T_h$ of the domain $\Omega$ of characteristic size $h$, and let us denote with $V_h$ a finite element space defined on $\mathcal T_h$. 
The dimension of the discrete space $V_h$ is denoted by $\cN$, i.e., the number of degrees of freedom per discrete (nodal) point of the finite element mesh.
Moreover, it is also assumed that $V_h$ satisfies suitable approximation properties with respect to the ambient space $V$. 

In what follows, $(\cdot,\cdot)$ will denote the inner product for the particular case of $L^2(\Omega)$.
Moreover, $\Pi_U$ will denote the orthogonal projector on a subspace $U \subset V_h$.
Finally, $\cL_i$, $i=1\hdots,\cN$ will denote the basis functions of the finite element space, while $M \in \bR^{\cN \times \cN}$ will stand for the mass matrix that defines the Hilbertian setting of the considered finite element space and the corresponding norm, i.e., for the $L^2$-setting:

\begin{equation}
M_{ij}  \coloneqq \int \cL_i, \cL_j \dx
\,.
\label{eq:M-matrix}
\end{equation}

\subsection{The space of observations}\label{ssec:W-space}

The first step to formalize the data assimilation problem consists of constructing 
a functional space of observations.
To this purpose, let  $l_1\hdots,l_m$ be a set of $m$ independent observations, assumed to be modeled by $m$ independent linear functionals $\ell_i: V_h \to \mathbb R$ $(i=1,\hdots, m)$ applied to an underlying \textit{true} solution $u_{\rm true} \in V_h$, i.e., $l_i = \ell_i(u_{\rm true})$, and let us gather them as entries of a vector of measures $l \in \bR^m$.

For instance, this might be the case of image voxels (or, more in general, any sensor data).

These assumptions allow to define a finite-dimensional space of observations $W = \text{Span}(w_i,\hdots,w_{m})$, with $W \subset V_h$, spanned by the unique Riesz representers of the observation functionals $\ell_i$,i.e., so that  it holds
\begin{equation}
\ell_i(u_{\text{true}}) = 
\< u_{\text{true}}, w_i \>, \; \mbox{for}\;i=1,\ldots,m
\,.
\label{eq:omega_i}
\end{equation}
\review{Since the measurements are assumed to be independent, 
without loss of generality, the Riesz representers are considered to be orthonormal, i.e., 
\begin{equation}
\< w_i, w_j \> = \delta_{ij}
\end{equation}
(orthonormality can otherwise be imposed in a post-processing step).}

Thus, the input data (i.e., the set of available measurements, $l_1,\hdots,l_m$) 
is modeled via a function $w \in W$ that is defined by the projection of the true solution onto $W$, i.e.,
$$w = \sum_{i=1}^m \ell_i(u_{\text{true}}) w_i = \Pi_{W} u_{\text{true}}.$$

\subsection{Data assimilation as a minimization problem}
\review{The data assimilation problem under investigation can be formalized as a linear reconstruction algorithm
for a \textit{measure-to-state} operator 
\begin{equation}\label{eq:Aop}
\cA : W \rightarrow V_h,
\end{equation}
which, given a set of observation $w \in W$ ($m$-dimensional space of observations), finds a suitable state 
$u^* \in V_h$ (ambient space, whose dimension depends on the finite element discretization) that approximates $u_{\text{true}}$, in the sense
that $\norm{u^* - u_{\text{true}}}$ is small.
}

The problem of estimating $u^* \in V_h$ solely relying on the data $w \in W$ is ill-posed. 
To overcome this issue, information about the dynamics of the state $u_{\text{true}}$ encoded in the forward problem \eqref{eq:model-pde} shall be incorporated.
Namely, a manifold of solutions (the `{parametrized-background}') of the underlying dynamics \eqref{eq:model-pde} is introduced

\begin{equation}
\cM = \{ u \in V_h;~\cP(u,\theta) = 0,~\theta \in \Theta \}
\,,
\label{eq:manifold_of_solutions}
\end{equation}
encompassing a range of governing dynamics as broad as the dimension of the parameter space $\Theta \subset \mathbb R^p$.
The dimension of the manifold $\cM$ depends on 
size of the discretization $\cN$ employed to numerically solve the underlying dynamics \review{and on the sampling used for the
parameter space}.

The parameter space $\Theta$ used to define the manifold \eqref{eq:manifold_of_solutions} can be 
chosen a subset of the admissible parameter space for the forward model (more details will be provided in 
Section \ref{ssec:training}). 
Therefore, one can at this stage incorporate in the data assimilation framework any additional knowledge on the 
problem parameters such as, e.g., initial estimates or variability relevant for the application of interest.

\review{The physical model is incorporated into the reconstruction problem via a reduced-order model.
Namely, under the assumption that the resulting dynamic has a fast decaying Kolmogorov $n$-width and, therefore, that the dynamics can be well approximated with a suitable reduced-order model, we introduce a hierarchy of (reduced-order) nested linear sub-spaces 
$V_n \subset \cM$, $n=1,2,\hdots$, with dimension $n \ll \text{dim}\,\cM$, that approximate well the manifold of solutions $\cM$.
Formally, we assume that the approximation errors of the reduced-order spaces
\begin{equation}\label{eq:rom_approximation}
\epsilon(V_n) :=\sup_{u \in \cM} \norm{\Pi_{V_n} u  - u}
\end{equation}
decay rapidly with $n$ (see Section \ref{ssec:training}).
}
\review{Hence, under the assumptions:
\begin{itemize}
\item The $\ell_i$ are linear functionals,
\item $V_n$ is a linear (or affine) sub-space of $V_h$,
\item The dimension of $V_n$ is smaller than the number of measures, i.e., $m > n$.
\item $V_h$ is a Hilbert space,
\end{itemize}
the reconstruction of the measure-to-state operator is formulated as a 
constrained minimization problem for the distance between the reconstruction and the reduced-order space \cite{MPPY2015}:
\begin{problem}[PBDW]\label{pb:pbdw}
For a given $w \in W$, find $u^* = \mathcal A(w) \in V_h$ such that:
\begin{equation}
u^* = \arginf_{u \in  V_h} \norm{\Pi_{V_n^\perp} u}^2
\,,\; \mbox{such that} \; 
\< w_i, u \> = \ell_i(u),\; i=1,\hdots,m
\,.
\label{eq:pbdw2}
\end{equation}
\end{problem}}

\review{The rest of this Section overviews the relevant theoretical and practical aspects concerning the solvability of problem 
\ref{pb:pbdw}. For a detailed analysis, including its convergence properties for increasing $m$ (i.e., the number of available observations), we refer to \cite{MPPY2015}.}

\subsubsection{Training phase}\label{ssec:training}
\review{The training phase encompasses the necessary steps to create a suitable solution manifold and
the hierarchy of reduced-spaces $(V_n)_{n}$, on which the problem \ref{pb:pbdw} can be formulated.

Let us consider $K$ samples from the parameter set, denoted as $y_1,\hdots,y_K \in \Theta$.
Also, let us introduce a set 
$$
\cM^{\text{training}} := \left\{ u^1, \ldots, u^K \mid \cP(u^i,y_i) = 0,\;i=1,\ldots,K \right\} \subset \cM,
$$ 
containing `snapshots' of the governing dynamics, i.e., 
finite element solutions of the underlying PDE for the considered parameter samples.

The reduced space $V_n \subset V_h$ is constructed by means of a principal component analysis (PCA; also called proper orthogonal decomposition, or POD, see \cite{rathinam2003new}). 
The basis of this space is computed by introducing the covariance matrix $C \in \bR^{K \times K}$ whose entries are given by $C_{ij} = \< u^i, u^j \>$, and  computing an eigenvalue decomposition 
\begin{equation}\label{eq:rom-eigen}
C = B \Lambda B^T
\end{equation}
for $\Lambda \in \bR^{K \times K}$ and $B \in \bR^{K \times K}$.

The reduced basis $\(\rho_1 | \ldots | \rho_n\) \in \bR^{\cN \times n}$
that spans $V_n$ can be then obtained from the first $n$ columns -- 
ordered according to the eigenvalues (from highest to lowest) -- of the matrix 
$$
U = S^{-1} A B,
$$ 
where $S^2 = \text{diag}\{\Lambda_i, i=1,\dots,\min\{K,\cN\}\}$, and $A$ is the so-called \textit{snapshots matrix}, 
whose columns correspond to the snapshots $u^i$, $i=1,\hdots,K$. }

\begin{remark}
The approach described above to compute the reduced basis is computationally convenient if $K \ll \cN$, which is the case for numerical experiments presented in the following sections.
When $K \gg \cN$ the PCA is typically performed with the covariance matrix $C = \sum_{i=1}^K u_i \otimes u_i$, which has dimension $\cN \times \cN$. 
\end{remark}

\review{
\subsubsection{Optimality conditions of the PBDW}\label{ssec:optimality}

Let $\Phi = \(\rho_1 | \ldots | \rho_n\) \in \bR^{\cN \times n}$ be the matrix whose columns contains the reduced basis. 
The projection operator $\Pi_{V_n^h}: V_h \to V_n$ can be written in matrix form as
\begin{equation}
\Pi_{V_n^h} u = \Phi \Phi^T M u, \; \forall{u \in V_h}
\,.
\label{eq:PiVn}
\end{equation}

Similarly, introducing the matrix $\cW = \(w_1 | \ldots | w_m\) \in \bR^{\cN \times m}$ whose columns are the Riesz representers of the measurements, the projection onto $W$ can be expressed via 
\begin{equation}
\Pi_{W} u = \cW \cW^T M u, \; \forall{u \in V_h}
\,.
\label{eq:PiWh}
\end{equation}
Notice that $\cW$ is orthogonal since orthonormality has been enforced on the representers $w_1,\hdots,w_m$.

Let $l \coloneqq (\ell_i(u_{\rm true}))_i \in \mathbb R^m$ be the vector of available measurements.
Using \eqref{eq:PiVn} and \eqref{eq:PiWh}, the Lagrangian for the minimization problem \eqref{eq:pbdw2} 
reads:
\begin{align}
\cL(u;\lambda)
&= \norm{u - \Pi_{V_n^h} u}^2 - \lambda^T(\cW^T M u - l)
\nonumber\\
&= \( \( I - \Phi \Phi^T M\) u \)^T M \( \( I - \Phi \Phi^T M\) u \) - \lambda^T\(\cW^T u - l\)
\nonumber\\
&= u^T M \( I - \Phi \Phi^T M\) u - \lambda^T \(\cW^T M u - l\)
\,,
\label{eq:pbdw_lagrangian}
\end{align}
with $\lambda \in \bR^m$. The corresponding Euler-Lagrange equations for \eqref{eq:pbdw_lagrangian} read
\begin{equation}
\begin{pmatrix}
I - \Phi\Phi^T M & -\cW \\
\cW^T            & 0
\end{pmatrix}
\begin{pmatrix}
u \\
\lambda
\end{pmatrix} 
=
\begin{pmatrix}
0 \\
l
\end{pmatrix}
\,,
\label{eq:EL}
\end{equation}
with $\lambda \in \bR^m$. 

The optimality condition for problem \eqref{eq:pbdw2} is hence expressed by the
saddle-point problem \eqref{eq:EL}, which has a unique solution if $n\leq m$ and if the following condition is satisfied

\begin{equation}
\beta(V_n,W) := \inf_{v \in V_n} \frac{\norm{\Pi_W v}}{\norm{v}} > 0
\,.
\label{eq:pbdw_bound}
\end{equation}
}
\begin{remark}\label{rem:beta}
\rrreview{The stability constant $\beta(V_n,W)$  depends only on the mutual structures of the 
reduced-order space $V_n$ ($n$-dimensional, spanned by $\rho_1,\hdots,\rho_n$) and of the space of observations $W$ 
($m$-dimensional, spanned by the functions $w_1,\hdots,w_m$), and, geometrically, 
it can be seen as the \textit{angle} between $V_n$ and $W$.
In particular, notice that  we can compute $\beta(V_n,W)$ by evaluating the minimal singular value of the matrix $G^T G$, 
where 
\begin{equation}\label{eq:G}
G = \cW^T M \Phi \in \bR^{m \times n}, \mbox{i.e.,}\; G_{ij} = \< w_i, \rho_j \>,
\end{equation}
therefore allowing us to verify the assumption \eqref{eq:pbdw_bound}, thus the well-posedness of the problem, beforehand.}
\end{remark}

Otherwise stated, $\beta(V_n,W)$ measures the \textit{observability} of the reduced space $V_n$ based on the type of available measurements.
Namely, for a fixed number of observations, $m$, the richer the reduced space is (i.e., the larger its dimension $n$), the more likely it will be to find elements that are hardly observable with the given set of  measurement -- or, equivalently, which lie orthogonal to $W$ -- hence deteriorating the algorithm stability according to \eqref{eq:pbdw-error-bound}.

An alternative to solving the saddle-point problem \eqref{eq:EL} for $(u,\lambda)$ is stated in the following result.
\review{\begin{proposition}\label{prop:uveta}
Assume that the condition \eqref{eq:pbdw_bound} holds. 
Then, the solution to the problem \eqref{eq:pbdw2} admits a decomposition of the form 
\begin{equation}\label{eq:uveta}
u^* = v^* + \eta^*,
\end{equation} 
with
\begin{equation}
v^* = \Phi \[ \( G^T G \)^{-1} G^T l \] \in V_n
\,,
\label{eq:v_star}
\end{equation}
where \rrreview{$G$ is the matrix defined in \eqref{eq:G}}
and
\begin{equation}
\eta^* = \cW (\cW^T M v^* - l) \in W
\,.
\label{eq:eta_star}
\end{equation}
Moreover, \eqref{eq:uveta} is an orthogonal decomposition, i.e., $\eta \in V_n^{\perp}$.

\end{proposition}
}
\begin{proof}
Let $(u,\lambda)$ be a solution of \eqref{eq:EL}. We consider
the orthogonal decomposition $u = v + \eta$, with  $v \in V_n$ and $\eta \in V^{\perp}_n = \text{Ker}(\Pi_{V_n})$. 
Since $\eta = \Pi_{V_n^{\perp}} u$, using \eqref{eq:EL} and the orthogonality of $\cW$, we obtain
\begin{equation}
\(I - \Phi \Phi^T M\)u - 
M \cW \lambda =
\eta - \cW \lambda = 
0 \Rightarrow 
\cW^T \eta =
\lambda
\,.
\label{eq:lam1}
\end{equation}
Multiplying by $\Phi^T M$ leads to:
\begin{equation}
0 = \Phi^T M \left(\eta - \cW \lambda\right) = 
\Phi^T M \cW \lambda = G^T  \lambda
\,,
\label{eq:lam2}
\end{equation}
where we used the definition $G = \cW^T M \Phi$.

Let now $c = (c_1,\hdots,c_n) \in \bR^n$ be the coordinates of $v$ in the basis $\Phi$, i.e. $v = \Phi\,c$.
Inserting the decomposition $u = \Phi\,c + \eta$ into \eqref{eq:EL} we obtain
\begin{equation}
\underbrace{\cW^T M \Phi\,}_{G} c + \cW^T \eta = l
\,.
\label{eq:proof-1}
\end{equation}
Multiplication by $G^T$ yields
\begin{equation}
G^T l =
G^T G \,c + G^T \cW^T \eta
\underbrace{=}_{\tiny \mbox{\eqref{eq:lam1}}}
G^T G \,c + G^T \lambda
\underbrace{=}_{\tiny \mbox{\eqref{eq:lam2}}}
G^T G \,c
\,.
\label{eq:lam3}
\end{equation}
Solving for $c$ in \eqref{eq:lam3} and using $v = \Phi\,c$ yields \eqref{eq:v_star}, while Equation \eqref{eq:eta_star} follows from \eqref{eq:proof-1}.
\end{proof}

\review{
\begin{remark}
Notice the presence of the pseudoinverse \rreview{$(G^T G)^{-1}G^T$} in equation \eqref{eq:v_star}. Hence, for a given $u \in V$, $v^*$ solves the least squares problem 
(see, e.g., \cite{MPPY2015}):
\begin{equation}
\inf_{v \in V_n} \frac12 \norm{\Pi_{W} v - w}^2,
\end{equation}
where $w = \Pi_{W} u$ denotes the vector of observations taken on $u$. 

\end{remark}
\begin{remark}
The sought operator $\cA(w)$ can be thus decomposed in a prediction, 
$v^* \in V_n$, and a correction, $\eta^* \in V_n^{\perp} \oplus W$, component respectively.
The former belongs to the considered reduced-order model of the parametric dynamics, while the latter depends on the discrepancy between the prediction and the data and can be seen as a corrector for any model bias present on $V_n$.
\end{remark}

From proposition \ref{prop:uveta},
it follows that, for the optimal state, it holds $u^* \in V_n \oplus \( W \cap V_n^\perp \)$, i.e., 
the measure-to-state operator $\cA(w)$ is a bounded linear map between $W$ and $V_n \oplus \( W \cap V_n^\perp \)$.
Using \eqref{eq:pbdw_bound} and \eqref{eq:rom_approximation}, 
one obtains that, for any $u \in V_h$, the reconstruction error can be bound by \cite{BCDDPW2017}
\begin{equation}\label{eq:pbdw-error-bound}
\begin{aligned}
\norm{u - \cA(w) } & \leq 
\beta(V_n,W)^{-1} \norm{u - \Pi_{V_n \oplus \( W \cap V_n^\perp \)}u } 
\leq \beta(V_n,W)^{-1}  \norm{u- \Pi_{V_n}u} \\
& \leq \beta(V_n,W)^{-1} \epsilon(V_n).
\end{aligned}
\end{equation}
}

In practice, denoting with $\Lambda_1 \geq \Lambda_2 \geq \ldots \geq \Lambda_{\max{\{K,\cN\}}}> 0$
the eigenvalues of the covariance matrix \eqref{eq:rom-eigen} (in decreasing order), the quality
error of the reduced-order model can be approximated as
\begin{equation}
\hat{\epsilon}_n = 
\(\sum_{i=1}^{\max\{\cN,K\}} \Lambda_i \)^{-1/2} \(\sum_{i=n+1}^{\max\{\cN,K\}} \Lambda_i \)^{1/2}\,.
\label{eq:mor_tails}
\end{equation}

Estimate \eqref{eq:pbdw-error-bound} indicates a trade-off between the approximation quality of the reduced-order model,
$\epsilon(V_n)$, and the observability (stability constant), $\beta(V_n, W)$, since both
quantities decrease with $n$.
This relation can be used to choose the dimension of the space $V_n$. 
For instance, in \cite{BCDDPW2017}, a so-called \textsl{nested space} strategy is followed to select the dimension of the reduced-order model so that the \textit{a priori} bound \eqref{eq:pbdw-error-bound} is minimized.

\rreview{
\begin{remark}
The results described in Section \eqref{ssec:optimality} are derived regardless from the physical model (i.e., the PDE) under consideration
and rely only on the Hilbert structure of the ambient space and on the hypothesis that the observations can
be represented as independent linear functionals of the state variables \cite{maday-pbdw-2014,GLM2021}. In particular, 
the theory covers also the cases of a partially observable state, i.e., when one component of the state variable cannot be observed.
One application of the framework in this situation will be discussed in more details in Section \ref{sec:mre-assimilation}, investigating the joint
reconstruction of displacement and pressure fields using displacement dependent observations only.
\end{remark}
}

\section{Poroelastic mechanics}\label{sec:forward}

This paper targets the application of the data assimilation framework described in Section \ref{sec:pbdw} to the case of brain elastography.
To this purpose, we model the brain tissue as a poroelastic media, composed of a biphasic mixture of solid gray matter, solid white matter,
 and cerebro-spinal fluid (CSF). 
Brain biomechanics is assumed to obey Biot's equations of poroelasticity \cite{perrinez_etal_2009}, the theory of which describes the behavior of the tissue at the mesoscale considering the interaction between biphasic material strains and increments in fluid volume. 

\subsection{Brain poroelastic model}\label{ssec:poroelasticity}

\begin{figure}[!htbp]
\centering
\includegraphics[height=5cm]{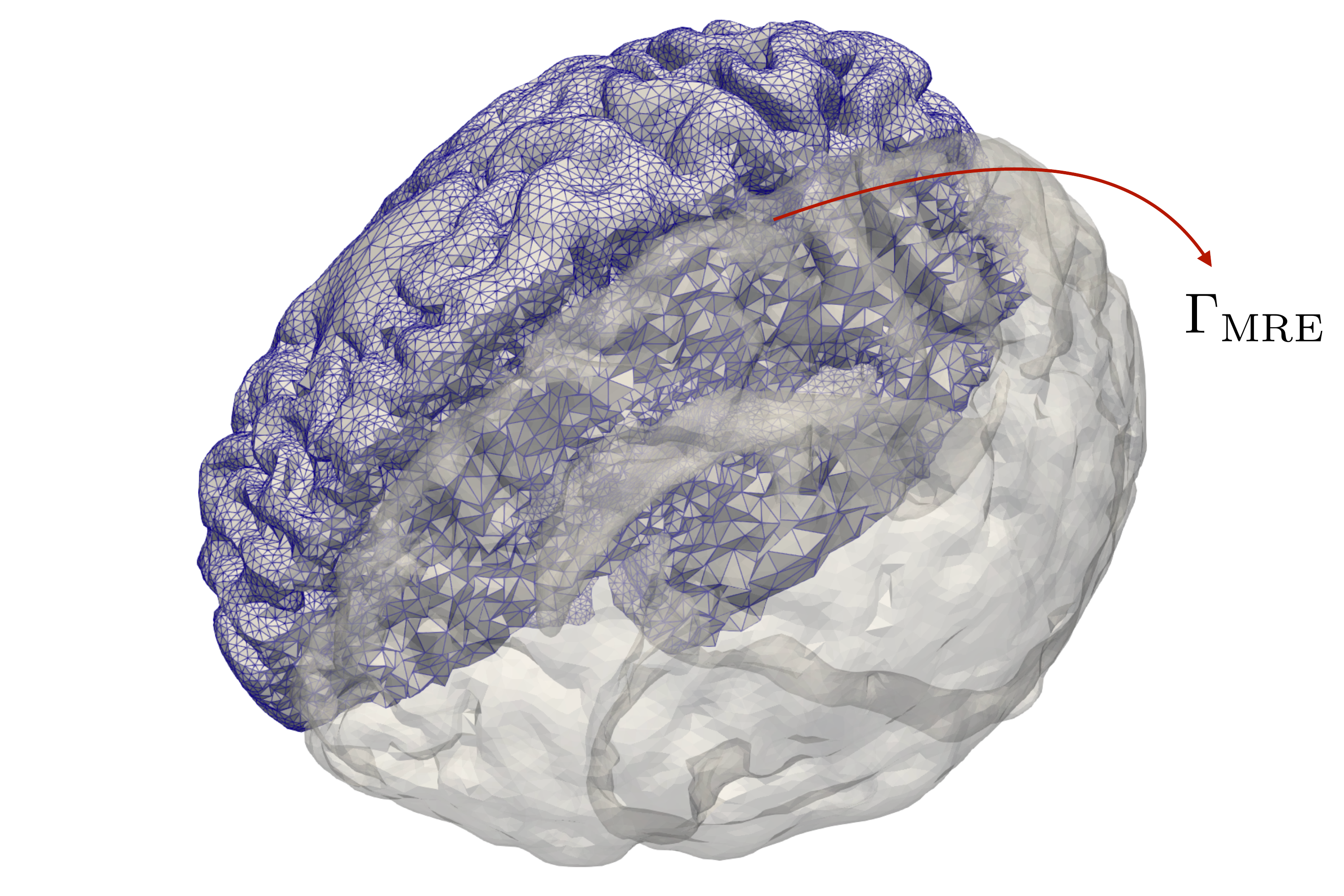}
\includegraphics[height=5cm]{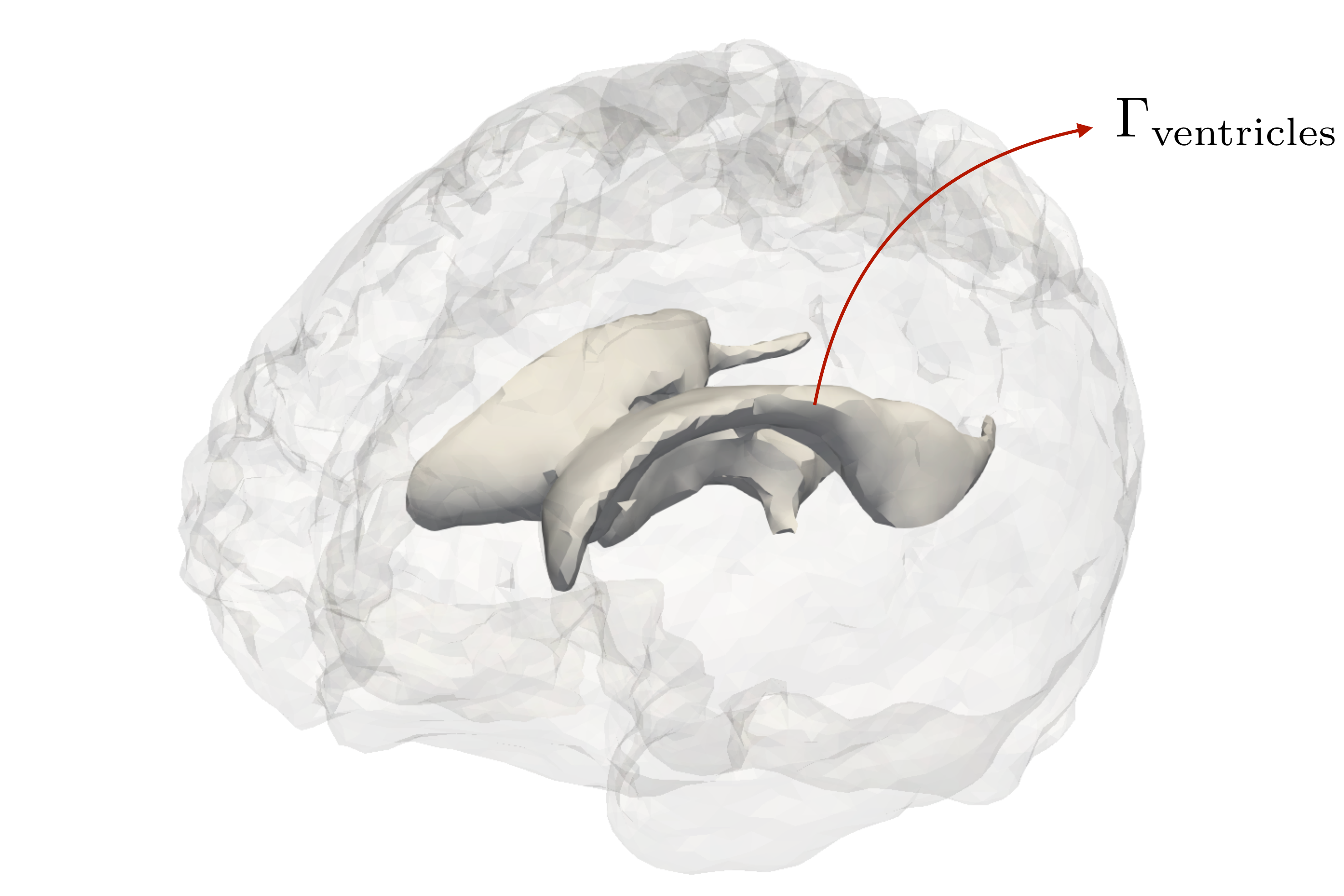}
\includegraphics[height=5cm]{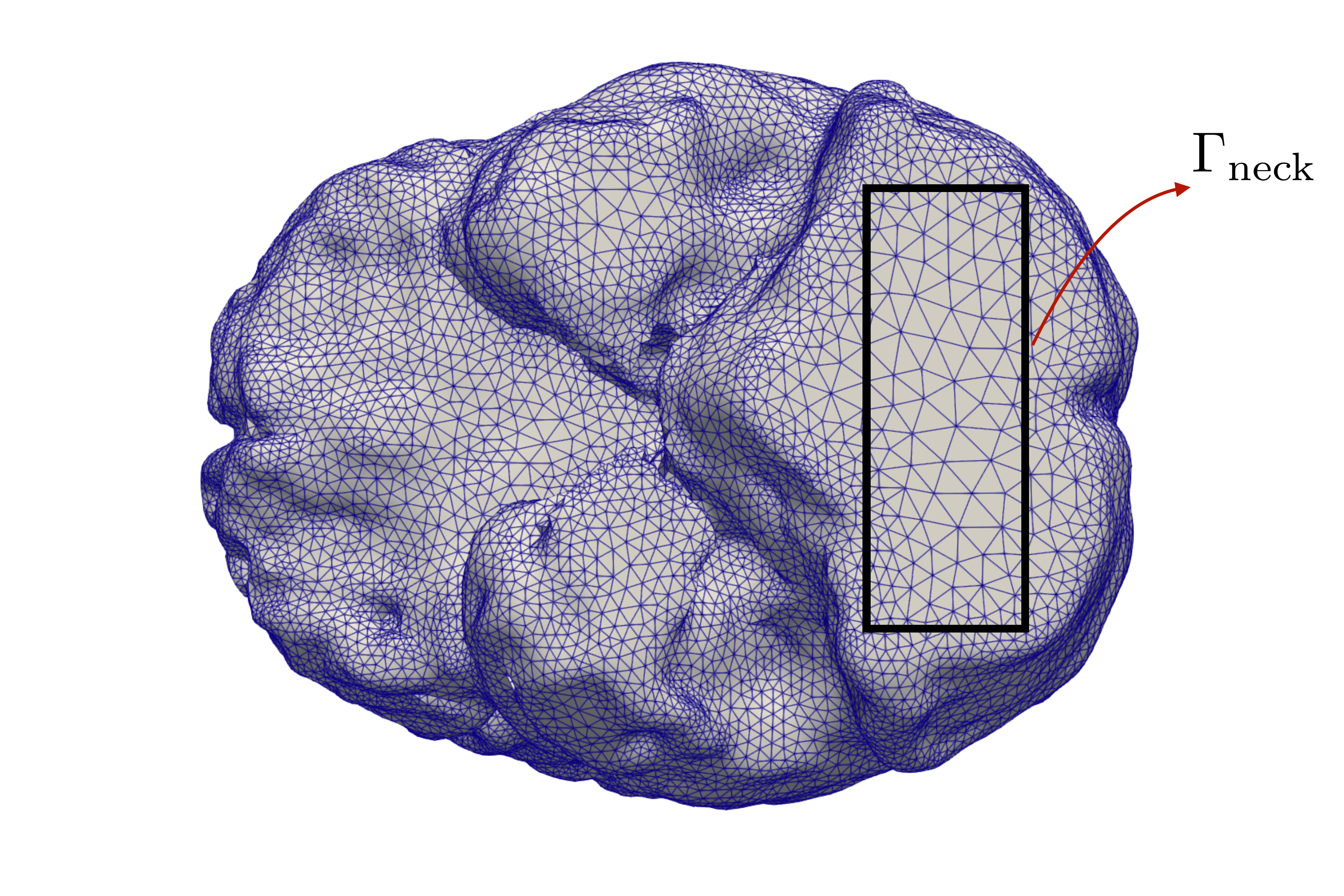}
\caption{Working domain $\Omega$ and the boundary decomposition.}
\label{fig:brain-forward-geo}
\end{figure}

Derivation of the poroelasticity equations is briefly outlined in this paragraph. 
Let $\Omega \subset \bR^3$ denote the computational domain, whose boundary $\partial \Omega$ is decomposed into disjoint sets: $\partial \Omega = \Gamma_{{\rm neck}} \cup \Gamma_{\text{ventricles}} \cup \Gamma_{{\rm MRE}}$.
$\Gamma_{{\rm neck}}$ denotes the portion of the boundary where the displacement vector field is zero, $\Gamma_{{\rm MRE}}$ is the boundary where a harmonic pulse (to replicate an MRE pulse) is prescribed, and $\Gamma_{\text{ventricles}}$ stands for the internal boundary between brain tissue and ventricles (see Figure \ref{fig:brain-forward-geo}).

Let $T$ denote the length of the considered time interval.
Soft tissue biomechanics is described by the solid displacements $\bu : \Omega \times [0,T] \rightarrow \bR^3$, the interstitial fluid pressure $p : \Omega \times [0,T] \rightarrow \bR$, and the filtration velocity $\bz : \Omega \times [0,T] \rightarrow \bR^3$ fields.
The governing equations of poroelasticity are derived from conservation laws for mass and momentum, considering the interplay between (solid and fluid) material strains and local increments in fluid volume.
The balance equations of linear momentum read

\begin{equation}
\rho\,\partial_{tt} \bu - \boldsymbol\nabla\!\cdot\!\boldsymbol\sigma(\bu)  + \alpha\,\boldsymbol\nabla p = \mathbf{0}
\quad \text{in } \Omega \times [0,T]
\,.
\label{eq:PEDup_momentum}
\end{equation}

\noindent
In \eqref{eq:PEDup_momentum}, $\rho$ stands for the solid mass density, and the stress (second-order symmetric) tensor $\boldsymbol\sigma(\bu)$ is related to the corresponding strain tensor by Hooke's law

\begin{equation}
\boldsymbol\sigma (\bu) =
\frac{E}{1+\nu}\,\boldsymbol\varepsilon +
\frac{E\nu}{(1+\nu)(1-2\nu)}\,{\rm tr}(\boldsymbol\varepsilon)\,\bI
\,,
\label{eq:sigma_ped_up}
\end{equation}
where $E$ and $\nu$ are Young's modulus and Poisson ratio of the tissue respectively.
Moreover, $\bI$ is the identity matrix in $\bR^d$, and $\boldsymbol\varepsilon$ denotes the symmetric strain deformation tensor, i.e., $\boldsymbol\varepsilon(\bu) = \frac12 \left( \boldsymbol\nabla \bu + \left(\boldsymbol\nabla \bu\right)^T\right)$.
Parameter $\alpha$ in \eqref{eq:PEDup_momentum} is the Biot-Willis parameter which quantifies the coupling between the stress due to the increment of interstitial pressure and the wave propagation stress in the solid matrix.
In Equation \eqref{eq:PEDup_momentum} we neglected the contribution of the fluid inertia, assuming that during MRE the relative displacement of the fluid is negligible.

Fluid flow in the porous medium is modeled using Darcy's law
\begin{equation}
\bz + \frac{\kappa}{\mu}\,\boldsymbol\nabla p = \boldsymbol{0}
\quad \text{in } \Omega \times [0,T]
\,,
\label{eq:PEDup_darcy}
\end{equation}
where $\kappa$ is the mixture permeability and $\mu$ is the fluid viscosity. 
The system of PDEs \eqref{eq:PEDup_momentum}--\eqref{eq:PEDup_darcy} is closed by a mass balance equation that accounts for the increment in fluid content and for fluid- and solid-phase velocities:

\begin{equation}
S_\epsilon\, \partial_t p + 
\boldsymbol\nabla\!\cdot \partial_t \bu + 
\boldsymbol\nabla\!\cdot \bz = 0
\quad \text{in } \Omega \times [0,T]
\,.
\label{eq:PEDup_mass}
\end{equation}

The mass storage parameter is described through: $S_\epsilon = {3\alpha(1-\alpha B)(1-2\nu)}{(B E)^{-1}}$, as a function of the Biot-Willis parameter
\review{$\alpha$}, the elastic parameters, $E$ and $\nu$, and the Skempton's parameter $B$.
For brain tissue described as a saturated medium, we set $\alpha = 1.$ and $B = 0.99$ \cite{smillie2004}. 

Taking the divergence of \eqref{eq:PEDup_darcy} and inserting the result into \eqref{eq:PEDup_mass} leads to the elimination of the fluid velocity $\bz$.
Together with \eqref{eq:PEDup_momentum}, one obtains the following system of PDEs  for the displacement and the pressure fields:

\begin{equation}
\begin{aligned}
\rho\, \partial_{tt} \bu - 
\boldsymbol\nabla\!\cdot \boldsymbol\sigma(\bu) +
\alpha\, \boldsymbol\nabla p = 0 &
\quad \text{in } \Omega \times [0,T]
\\
S_\epsilon\, \partial_t p +
\alpha\, \boldsymbol\nabla\!\cdot \partial_t \bu - 
\frac{\kappa}{\mu}\,\nabla^2 p = 0 &
\quad \text{in } \Omega \times [0,T]
\end{aligned}
\,.
\label{eq:the_model_PEDup}
\end{equation}

\begin{figure}[!htbp]
\centering
\includegraphics[height=7cm]{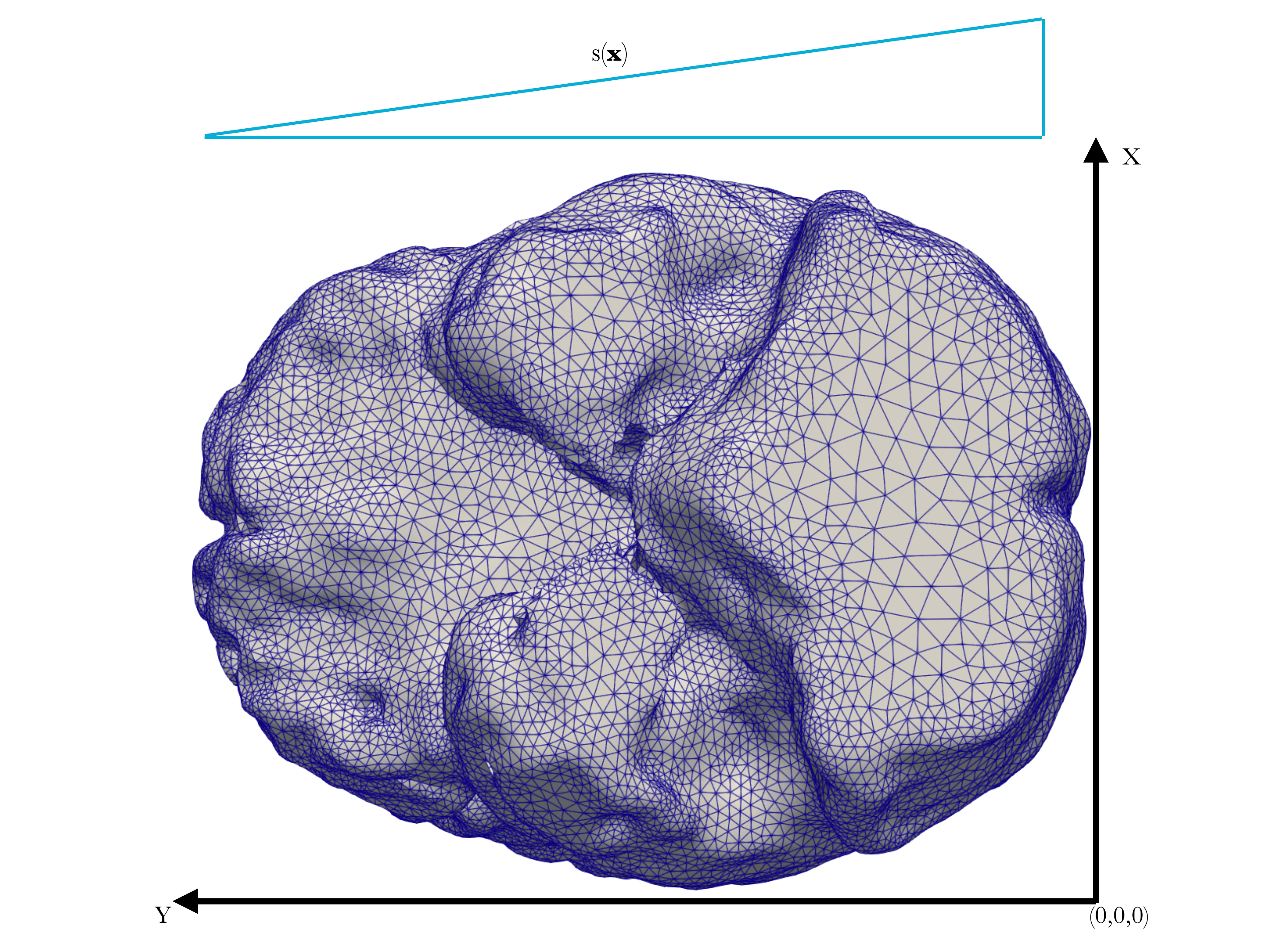}
\caption{Sketch of the model used for boundary conditions. To mimic the MRE setting, a harmonic force with a given amplitude and frequency has been multiplied by a linearly decreasing function $s(x)$ whose values (blue sketch on top) range from one -- where the mechanical pulse is applied -- to zero, along the direction of the pulse force.}
\label{fig:sketch_g}
\end{figure}

The system is closed by homogeneous initial conditions and by the following boundary conditions
\begin{equation}
\begin{aligned}
\boldsymbol\sigma \cdot \bn = \bg_{\rm MRE} &
\quad \text{on } \Gamma_{\text{MRE}}
\\
p = p_{\text{ventricles}} &
\quad \text{on } \Gamma_{\text{ventricles}}
\\
\bu = \mathbf 0 &
\quad \text{on } \Gamma_{\text{neck}}
\\
p = p_{\text{csf}} &
\quad \text{on } \Gamma_{\text{MRE}} 
\end{aligned}
\,,
\label{eq:the_model_PEDup_bc}
\end{equation}
with $p_{\text{csf}}$ and $p_{\text{ventricles}}$ denoting the CSF pressure imposed on the outer CSF and on the ventricles, respectively, while the load $\bg_{\rm MRE}$ imposed for a Neumann boundary condition on $\Gamma_{\text{MRE}}$ has been defined to mimic the forces on the brain tissue in an
 MRE examination.
Thus, the forcing term has been modeled via: $\bg_{\text{MRE}}(\boldsymbol x,t) = \xi \sin \(2 \pi w t\) s(\boldsymbol x)$, i.e., a mono-harmonic pulse (with a given cyclic frequency $\omega$ and amplitude $\xi$) multiplied by a space-depending function $s(\boldsymbol x)$ equal to unity where the pulse is applied (back of the brain) and linearly decreasing towards the front of the brain (see Figure \ref{fig:sketch_g}), i.e., $s(x,y,z) = \(1 - \nicefrac{y}{L}\)$, where $L=16.86$ cm is the geometry length along the y-axis (frontal axis).
This model is motivated by the need to account for the effect of the rigid skull on the propagation of the mechanical force over the whole brain surface.

The model parameters applicable for all simulations setup are summarized in Table \ref{tab:parameters}.

\begin{table}[!htbp]
\centering
\begin{tabular}{|c|c|c|c|}
\hline
Parameter & Description & Value [Unit] & Reference \\ \hline
$\alpha$ & Biot-Willis parameter & $1.0$ & \cite{smillie2004} \\ 
$\mu$ & Fluid viscosity & $10^{-2}$ $P$ & \cite{smillie2004} \\
$\rho$ & Solid matrix density & $1.0$ gr/cm\textsuperscript{3} & \cite{smillie2004} \\
$\kappa$ & Tissue permeability & $[10^{-9},~10^{-8}]$ cm\textsuperscript{2} & \cite{smillie2004} \\
$E$ & Solid matrix Young Modulus & $[10^5,~10^6]$ dyn/cm\textsuperscript{2} & \cite{smillie2004} \\
$\nu$ & Solid matrix Poisson's ratio & $[0.40,~0.45]$ & \cite{smillie2004} \\
$w_{\text{MRE}}$ & MRE frequency & $50$ Hz & \cite{hirsch2017magnetic} \\
$p_{\text{csf}}$ & CSF pressure & $10^4$ dyn/cm\textsuperscript{2} & \cite{bogo2016} \\
$p_{\text{ventricles}}$ & Ventricular pressure & $[1.0,~1.1]\times 10^4$ dyn/cm\textsuperscript{2} & \cite{bogo2016} \\
B & Skemptom parameter& 0.99 & \cite{smillie2004} \\
$\xi$ & Pulse magnitude & $-500$ dyn & $*$ \\ \hline
\end{tabular}
\caption{Parameters used in the poroelastic model and manifold generation. \\ 
($*$) The value of the pulse excitation magnitude, $\xi$, has been tuned so that brain tissue deformation simulation outputs are physiologically realistic.}
\label{tab:parameters}
\end{table}

\subsection{FE discretization and numerical solution}\label{ssec:fem} 

The system of equations \eqref{eq:the_model_PEDup} is numerically solved using the finite element method (FEM).
To this purpose, let $\mathcal T_h$ denote an unstructured tetrahedral  mesh that discretizes the computational domain $\Omega$ -- the characteristic size of the mesh is expressed through $h$.
We introduce the functional spaces 

\begin{equation}
\begin{aligned}
U & =  H^1(\Omega)^d_0 \,\coloneqq
\left\{ \bv \in H^1(\Omega)\mid \bv = 0 \text{ on } \Gamma_{\text{neck}} \right\}
\\
P & \coloneqq H^1(\Omega), ~~ Q \,\coloneqq 
\left\{ q \in H^1(\Omega) \mid q =  0 \;\mbox{on} \;\Gamma_{\text{ventricles}} \cup \Gamma_{\text{MRE}}\hdots \right\} 
\end{aligned}
\end{equation}

\noindent
and their discrete approximation by linear FE spaces 

\begin{equation}\label{eq:fe-spaces}
\begin{aligned}
U_h & \coloneqq 
\left\{ \bv \in U \mid \bv|_{T} \in \mathbb P_1(T), \forall T \in \mathcal T_h \right\}
\,, \\
P_h & \coloneqq 
\left\{ q \in P \mid q|_{T} \in \mathbb P_1(T), \forall T \in \mathcal T_h \right\}
\,, \\
Q_h & \coloneqq 
\left\{ q \in Q \mid q|_{T} \in \mathbb P_1(T), \forall T \in \mathcal T_h \right\}
\,.
\end{aligned}
\end{equation}

Then, we consider the following discrete weak formulation of \eqref{eq:the_model_PEDup}:
Find $(\bu_h, p_h) \in V_h \coloneqq U_h \times P_h$ such that:

\begin{equation}\label{eq:biot_semi_discrete}
\begin{aligned}
\rho \( \partial_{tt} \bu_h , \bv_h \) + 2E \( \boldsymbol\varepsilon (\bu_h), \boldsymbol\varepsilon(\bv_h) \) + \lambda \( \boldsymbol\nabla\!\cdot \bu_h,\boldsymbol\nabla\!\cdot \bv_h \)+ \alpha \( \boldsymbol\nabla p_h , \bv_h \) &= \( \bg_{\rm MRE}, \bv \)_{\Gamma_\text{MRE}}
\\
\( S_\epsilon\, \partial_t p_h, q_h \) + \(\alpha \boldsymbol\nabla \!\cdot \partial_t \bu_h, q_h \) +  \( \frac{\kappa}{\mu} \boldsymbol\nabla p_h, \boldsymbol\nabla q_h \) &= 0,
\end{aligned}
\end{equation}

\noindent
for all $\(\bv_h, p_h\) \in U_h \times Q_h$.

\review{
System  \eqref{eq:biot_semi_discrete} can be written in the following matrix form
\begin{equation}\label{eq:full-up-system}
\left[ 
\begin{array}{cc}
M_U & 0 \\ 0 & 0 
\end{array}
\right]
\left[ 
\begin{array}{c}
\partial_{tt}{\bu_h} \\ \partial_{tt}{p_h}
\end{array}
\right] +
\left[ 
\begin{array}{cc}
0 & 0 \\ \alpha B' & S_{\epsilon} M_P 
\end{array}
\right]
\left[ 
\begin{array}{c}
\partial_t{\bu_h} \\ \partial_t{p_h}
\end{array}
\right]+
\left[ 
\begin{array}{cc}
A_U & \alpha B  \\ 0 & \kappa \mu^{-1} A_P
\end{array}
\right]
\left[ 
\begin{array}{c}
{\bu_h} \\ {p_h}
\end{array}
\right] = 
\left[ 
\begin{array}{c}
F_U \\ F_P
\end{array}
\right]
\end{equation}
Let us introduce $y = [\bu_h, p_h]^T$,  $\dot{y} = [\partial_t \bu_h, \partial_t p_h]^T$, and 
$\ddot{y} = [\partial_{tt} \bu_h, \partial_{tt} p_h]^T$
and rewrite \eqref{eq:full-up-system} at time step $t_{n+1}$ as
\begin{equation}\label{eq:full-system-y}
M \ddot{y}_{n+1} + C \dot{y}_{n+1} + A y_{n+1} = F_{n+1}\,.
\end{equation}
}

\review{We used a Newmark method \cite{newmark} for the 
discretization of the second-order time derivatives, whereas a backward Euler algorithm is used for the first-order time derivatives.
Namely, for a fixed $0 \leq \hat{\beta} \leq \frac12$, we define
\begin{equation}\label{eq:newmark-y-dot-ddot}
\begin{aligned}
\ddot{y}_{n+1} &= \frac{1}{2\hat{\beta}} \frac{2 (y_{n+1} - y_n - \tau \dot{y}_n)}{\tau^2}
-\frac{1- 2\hat{\beta}}{2 \hat{\beta}} \ddot{y}_n \\
\dot{y}_{n+1} & = \dot{y}_n + \tau \ddot{y}_{n+1}
\end{aligned}
\end{equation}
where $\tau$ denotes the time step.
Inserting \eqref{eq:newmark-y-dot-ddot} into \eqref{eq:full-system-y}, the time-discrete formulation at time $n+1$ 
reads
\begin{equation}\label{eq:newmark-time-disc-y}
\begin{aligned}
&\frac{1}{\hat{\beta} \tau^2} M y_{n+1} + \frac{1}{\hat{\beta} \tau} C y_{n+1} + A y_{n+1} = F_{n+1} + 
\underbrace{\frac{1- 2\hat{\beta}}{2 \hat{\beta}} (M+C) \ddot{y}_n  +
 \left(\frac{1}{\hat{\beta} \tau^2}M +\frac{1}{\hat{\beta}\tau}C\right)(y_n + \tau \dot{y}_n)}_{:= G(y_n,\dot{y}_n,\ddot{y}_n)}
\end{aligned}
\end{equation}
where $G(y_n,\dot{y}_n,\ddot{y}_n)$ has been introduced to group all the terms dependent on the previous time iteration.}

\review{Multiplying by $\hat{\beta} \tau$ the second row of equation \eqref{eq:newmark-time-disc-y}, 
the corresponding problem in variational form reads: Find $(\bu_h^{n+1}, p_h^{n+1}) \in V_h \coloneqq U_h \times P_h$ such that
\begin{equation}\label{eq:forward-weak}
A( (\bu_h^{n+1},p_h^{n+1}),(\bv,q) )  = \( \bg_{\rm MRE}, \bv \)_{\Gamma_\text{MRE}} + 
\mathbf g_{\hat{\beta}}^n(\bv,q) \,,
\end{equation}
\noindent
for all $\(\bv_h, p_h\) \in U_h \times Q_h$, where the term $\mathbf g_{\hat{\beta}}^n$ depends on the previous iteration
$(\bu_h^n,p_h^n)$ and
\begin{align}
A( (\bu_h,p_h),(\bv_h,q_h) ) \coloneqq & \( \frac{\rho}{\hat{\beta}  \tau^2} \bu_h , \bv_h \) + 2E \( \boldsymbol\varepsilon (\bu_h), \epsilon(\bv_h) \) +
\lambda \( \boldsymbol\nabla\!\cdot \bu_h,\boldsymbol\nabla\!\cdot \bv_h \)
\nonumber \\
&~~~ + \alpha \( \boldsymbol\nabla p_h , \bv_h \) +  \( S_\epsilon p_h, q_h \) + \(\alpha \boldsymbol\nabla\!\cdot \bu_h, q_h \) + \hat{\beta} \tau\( \frac{\kappa}{\mu} \boldsymbol\nabla p_h, \boldsymbol\nabla q_h \)\label{eq:A-weak}\,.
\end{align}

Integrating by parts the term $\alpha \( \boldsymbol\nabla p_h, \bv_h \)$ testing the resulting formulation by $\( \bu_h, p_h \)$,
and using that $\hat{\beta} < 1$, one obtains the following stability results for the finite element formulation.}

\begin{proposition}
It holds 
\begin{equation*}
A( (\bu,p),(\bu,p) ) \geq 
\review{\hat{\beta} \|| (\bu,p) \||_{V_h}^2}
\,,
\end{equation*}
with the norm 
\begin{equation*}
\|| (\bv,q) \||_{V_h}^2 =
\frac{\rho}{\tau^2} \| \bv \|^2 + 
2E \| \boldsymbol\varepsilon(\bv) \|^2 + 
S_\epsilon \| p \|^2 + 
\frac{\kappa\,\tau}{\mu} \| \boldsymbol\nabla p \|^2
\,.
\end{equation*}
\end{proposition}

Two remarks are important.
First, it should be noted that the time derivative of the pressure and the diffusive term, $\( \nicefrac{\kappa\,\tau}{\mu} \boldsymbol\nabla p_h, \boldsymbol\nabla q_h \)$, act as a stabilizer of the spatial discretization.
Second, the resulting stability of $A( (\bu,p),(\bv,q) )$ strongly depends on the physical parameters, and it can indeed deteriorate -- especially concerning the pressure regularity -- if the mass storage coefficient $S$ or the medium permeability $\kappa$ are very small. 
This issue has also been investigated in \cite{rodrigo-2016} for the steady poroelastic case in the presence of strong discontinuities in the permeability. 
In that setting, it has been proposed to include a stabilization term of the form $\beta \( \boldsymbol\nabla p_h, \boldsymbol\nabla q_h \)$, where $\beta = \beta(E,\lambda,h)$.
In the physical regime relevant to this work, we did not observe stability issues related to spatial discretization. 
Additional numerical tests have been performed considering the stabilization proposed in \cite{rodrigo-2016} showing, however, only negligible differences.

\section{PDE-informed MRE data assimilation}\label{sec:mre-assimilation}

This section is devoted to the application of the data assimilation framework (outlined in Section \ref{sec:pbdw} from a general perspective) to the case of brain elastography.
In particular, it is presented how the displacement obtained via MRE will be assimilated into the poroelastic model detailed in Section \ref{sec:forward} and all the practical implementation aspects of the PBDW for the joint reconstruction of displacements and pressure fields are presented in detail. 
Further numerical results concerning the validation and application of the framework are presented in Section \ref{sec:results}.

The computations presented in these Sections have been performed using the software {MAD} \cite[Chapter 5]{galarceThesis}, which is based on the linear algebra library {PETSc} \cite{petsc}.

\subsection{Image acquisition and model generation}

The data assimilation pipeline is based on different sources of data:
\begin{itemize}
\item[(\emph{a})] The computational model has been generated using full brain anatomical, high-resolution MPRAGE MRI images (isotropic voxel size 1 mm\textsuperscript{3}).
The anatomical image data were acquired from a healthy volunteer after signing consent to use these for the present study.
\review{The acquisition time for the anatomical image was about 30 minutes}.
The raw DICOM data were then segmented using open-source software {3D Slicer} \cite{3dslicer} to produce a triangulated representation of the brain and ventricles surface (exported in an STL formatted file).
Subsequently, the surface mesh was loaded into {Mmg} \cite{mmg3d} to generate an unstructured mesh consisting of 12,345 4-node tetrahedral finite elements.
\item[(\emph{b})] The reconstruction is based on displacement data over a few slices of the three-dimensional brain, located as those acquired in an MRE examination. 
However, no assumptions on the spatial resolution of the displacement data, i.e., of the image voxels, are required.
\end{itemize}

The results shown in Section \ref{sec:results} are restricted to synthetic data, i.e., displacement fields generated from forward simulations on the physiological brain geometry. 
However, in the context of clinical routine, this type of data can be acquired via MRE.
During this examination the tissue is subject to a harmonic mechanical vibration (frequency 10-50 Hz) imposed by actuators, and the tissue response is recorded via phase-contrast MRI. The internal displacement field can then be extracted from the complex phase of the acquired images.

\review{For the case of synthetic measurements, it is not required to register the displacement data to the
patient geometry. However, in general, } 
the  definition of the space of observations (Section \ref{ssec:W-space}) requires the displacement data to be registered to the anatomical brain data. \review{In those situations, registration can be performed using validated software tools
for neuroimaging, see, e.g., \cite{GREVE200963}}.

Practically, the MPRAGE used for the definition of the computational mesh, and the MRE data to be assimilated, shall be acquired, on each patient, in a single examination, to minimize additional efforts for image registration.
This acquisition protocol has been recently used in the framework of inversion recovery MRE (IRMRE), to estimate biophysical parameters of the gray and white matter \cite{lilaj_etal_2021a}.
The application of the PBDW to in vivo data, especially within IRMRE, is the subject of ongoing work.

\subsection{Characterization of \review{observations}}
\begin{figure}[!htbp]
\centering
\subfigure[$x$-component (cm)]{
\includegraphics[height=4.2cm]{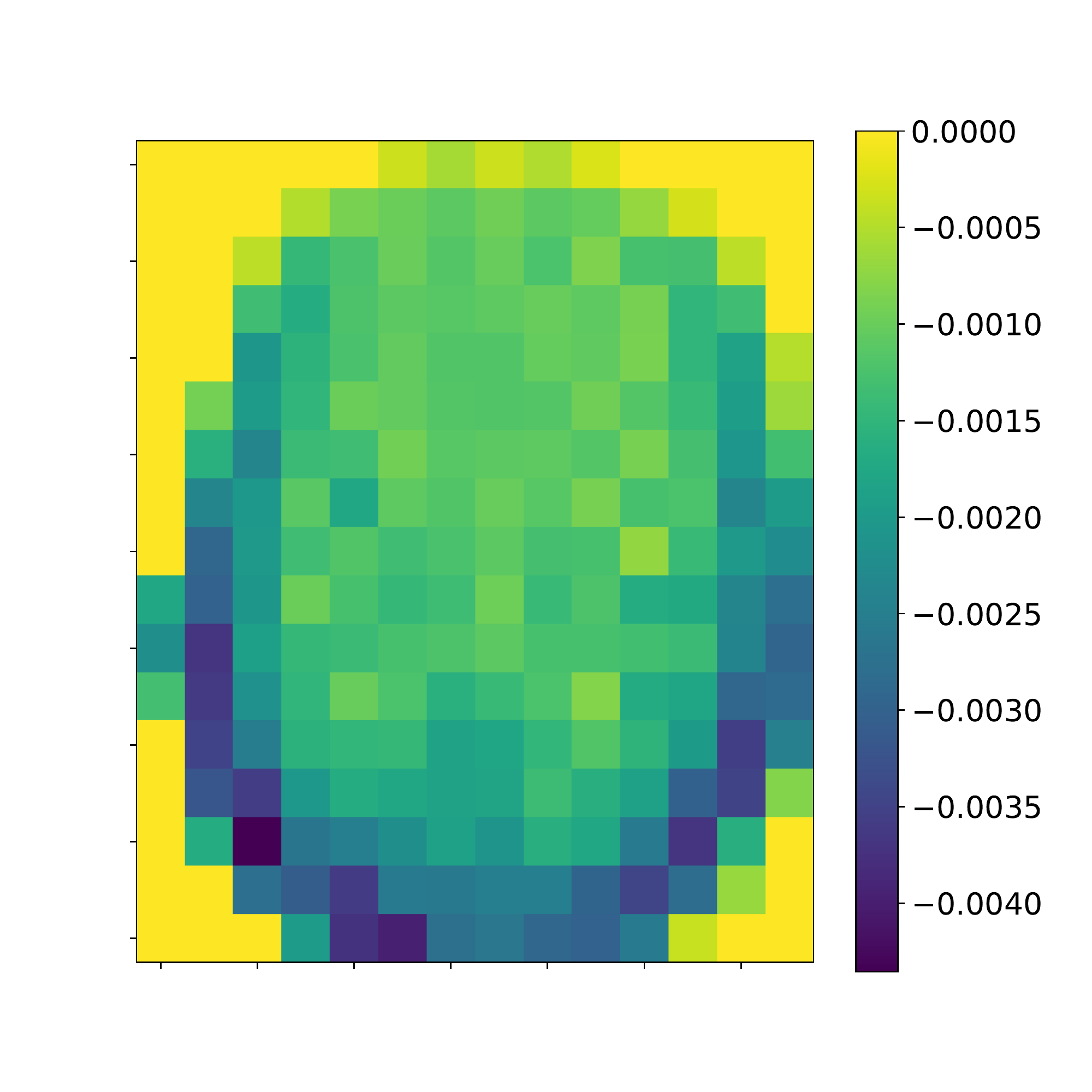}}
\subfigure[$y$-component (cm)]{
\includegraphics[height=4.2cm]{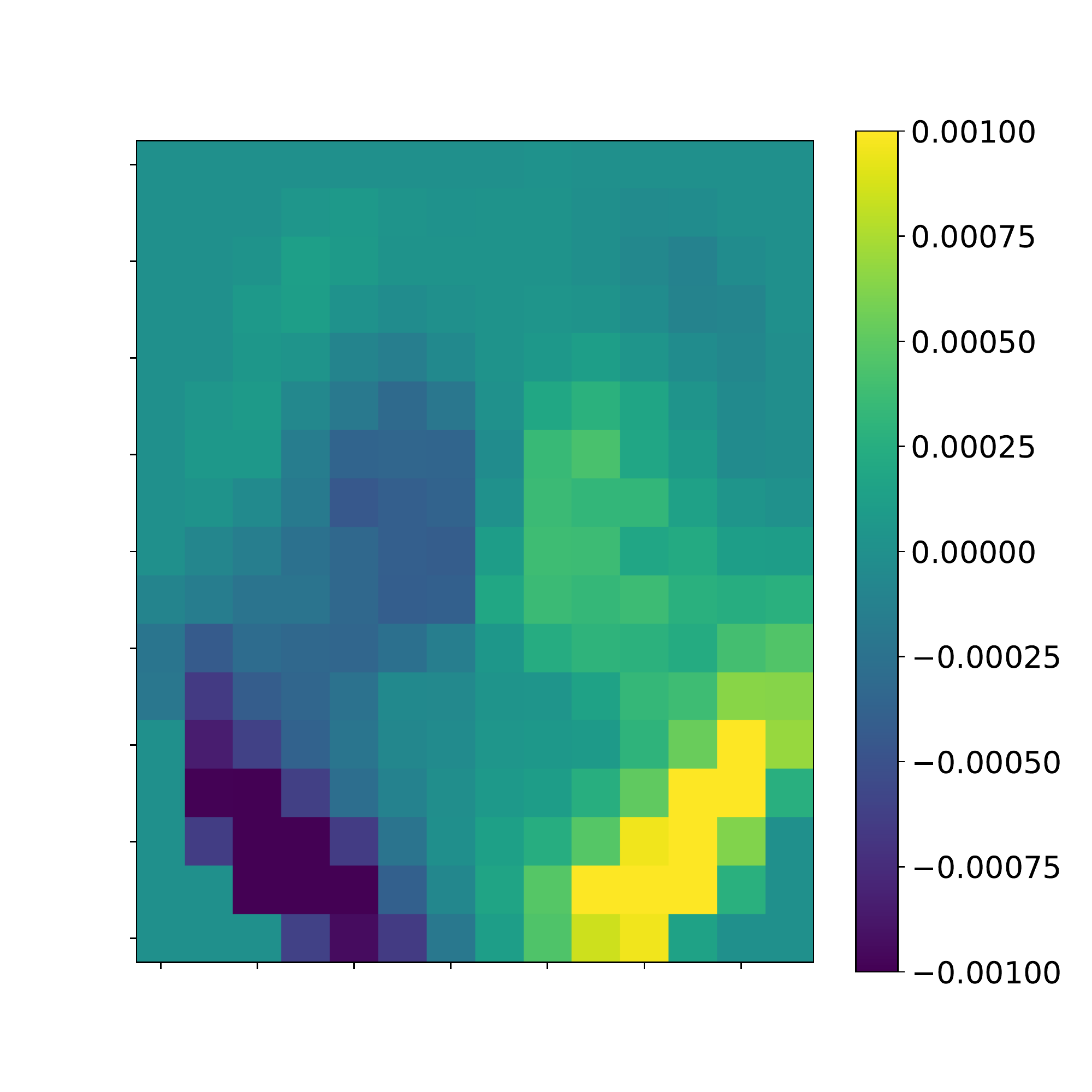}}
\subfigure[$z$-component (cm)]{
\includegraphics[height=4.2cm]{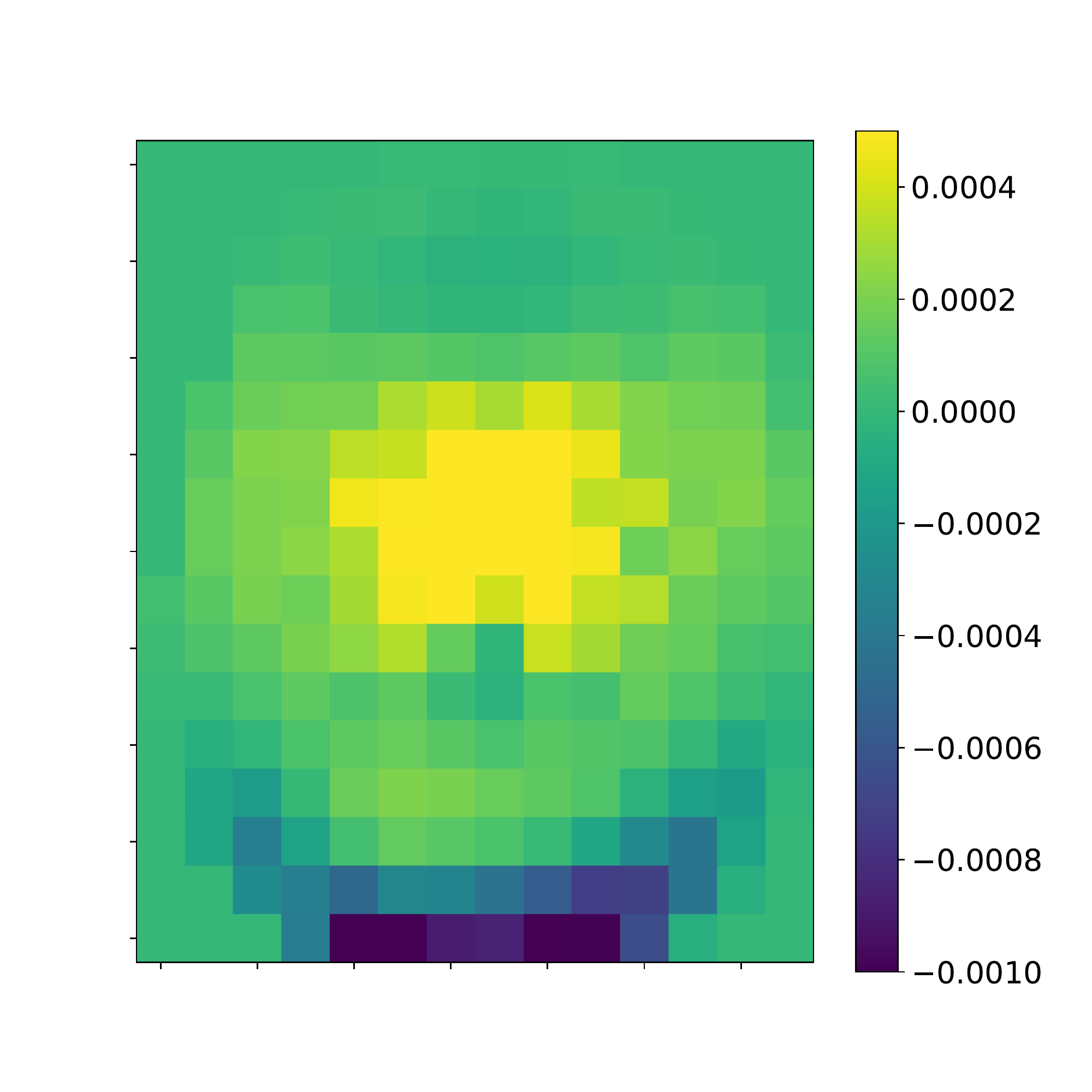}}
\subfigure[Location]{
\includegraphics[height=2.5cm]{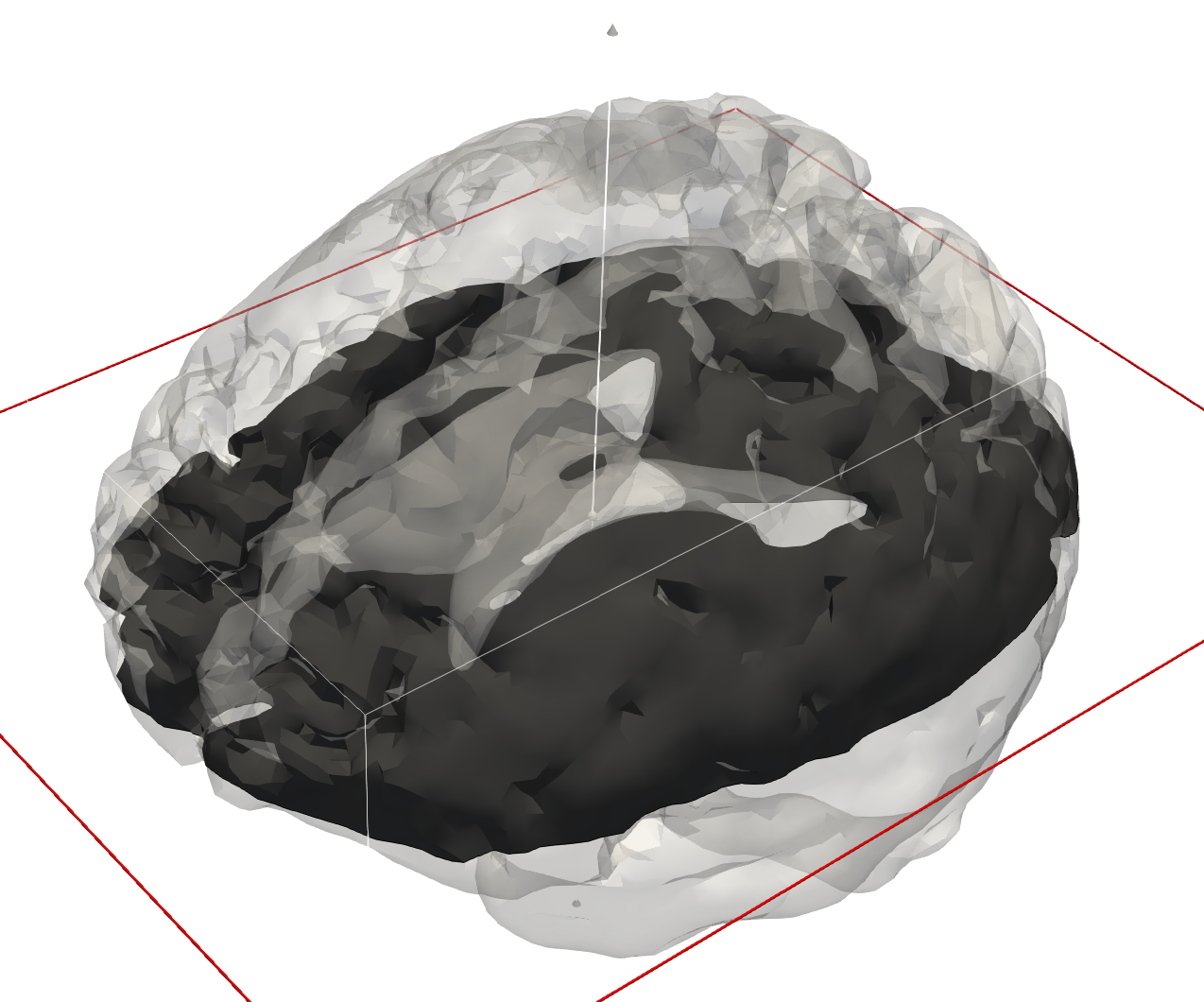}}
\caption{Figures (a), (b), and (c): \review{Selected} snapshot (time instant: 15 ms after the pulse begins) of the synthetic displacement images
considered for the data assimilation problem. Figure (d): Location of the measured slice with respect to the computational brain model.}
\label{fig:pbdw_measures}
\end{figure}

The available data are assumed to represent the three-dimensional displacement field on $N_v$ voxels in the upper part of the brain. Thus, each image contains $3\times N_v$ measurements (scalars), i.e., one 3D vector per imaged voxel (Figure \ref{fig:pbdw_measures}). 

In what follows, these voxels are denoted as $\Omega_i \subset \bR^3$ for $i = 1,\hdots,N_v$.
The $3 \times N_v$ linear functionals describing the observations are defined by taking the average component of the displacement vector, $v$, on each voxel. Namely, the basis of the space $W$ can be then computed after solving the following set of problems: Find a $w_i^j \in V_h$ such that

\begin{equation}
\< w_i^j , v \> = 
\ell_i^j (v) \coloneqq \int_{\Omega_i} v\,e_j\,\dx\,
\quad \forall{v \in V_h} ~
\mbox{for $i=1,\ldots,N_v$ and $j=1,2,3$}
\,,
\label{eq:rr_compute}
\end{equation}

\noindent
where $e_j$ is the $j$-th component of the unit vector.
\review{The linear systems \eqref{eq:rr_compute}, whose forms depend on the 
Hilbert structure that defines the scalar product on $V_h$, are solved using $\bP_1$ finite elements}.
The modified Gram-Schmidt algorithm is used for the ortho-normalization of the computed basis.

\begin{remark}[Offline computation of $W$]
\rreview{The space $W$ depends on the location of the imaged voxels, $\Omega_i$, but not on the particular value of the observations. 
Therefore, the computation of the basis of $W$, which requires the solution of $m$ variational problems  of size $\cN$
\eqref{eq:rr_compute},
does not depend on the particular patient under examination, but only on the setting of the image  acquisition device (i.e., on the location of the acquired data).
The space $W$ can be therefore computed \textit{offline}, i.e.\ before the remaining steps of the data assimilation algorithm (training, model reduction, and reconstruction) are performed. 
The \textit{online phase} of the algorithm only requires the solution of an $n \times n$ system of equations when the new observations are given.}
\end{remark}

A snapshot of the considered synthetic \review{measurements} is shown in Figure \ref{fig:pbdw_measures}.
The results shown in this work employ synthetic data, calculated computing a forward solution (Section \ref{ssec:fem}) for a selected set of parameters and evaluating local displacement averages on a voxel of 1$\times$1$\times$1 cm\textsuperscript{3} in size.
This procedure leads to an image containing 600 independent observations \review{on a brain slice}. 
\review{Further details on the observations used for the numerical experiments are provided
in  Section \ref{sec:results}.}

\begin{remark}[Resolution of MRE data]
The procedure for constructing the measurement space $W$ can be analogously defined for arbitrary image resolution, regardless of the  size of the finite element mesh elements.
\end{remark}

\begin{remark}
Using the formalism introduced in Section \ref{sec:pbdw}, for a given set of physical parameters these observations define the projection of the true state onto $W$, e.g.: $w = \Pi_W u_{\text{true}}$.
\end{remark}

\subsection{Training set and order reduction of the model}

The training manifold $\cM^{\text{training}}$ is computed sampling uniformly the parameter space in the intervals (see also Table \ref{tab:parameters}):
\begin{equation}
\kappa \in \[ 10^{-9}, 10^{-8} \]~\text{cm}^2
\,,\;
E \in \[ 10^5, 10^6 \]~\text{dyn}/\text{cm}^2
\,,\;
\nu \in \[ 0.4, 0.45 \]
\,,\;
p_{\text{ventricles}} \in \[ 1, 1.1 \]\times 10^4 dyn/\text{cm}^2
\,.
\label{eq:theta-training}
\end{equation}

The parameter ranges for $\kappa$, $E$, and $\nu$ used in \eqref{eq:theta-training} are based on  physiological parameter values obtained from relevant brain biomechanics modeling works from the literature.
All remaining parameters are assumed to be given and equal to the values specified in Table \ref{tab:parameters}.
For the easiness of notation, in what follows, the vector $\theta = \(\kappa, E, \nu, p_{\text{ventricles}}\)$ will denote a generic element of the parameter space.

The \rreview{CSF pressure is treated as an additional parameter and it} is introduced through the forward PDE model as a boundary condition. The considered range serves to be able to capture, within the training manifold, the variability of the ventricular CSF pressure between the two selected ranges (more details will be presented in Section \ref{ssec:pv-compute}).
The training manifold is computed numerically by solving together equations \eqref{eq:forward-weak}--\eqref{eq:A-weak} using the finite element method presented in Section \ref{ssec:fem}.  

\begin{figure}[!htbp]
\centering
\subfigure[Displacement field (magnitude) on the external surface]
    {\includegraphics[height=5cm]{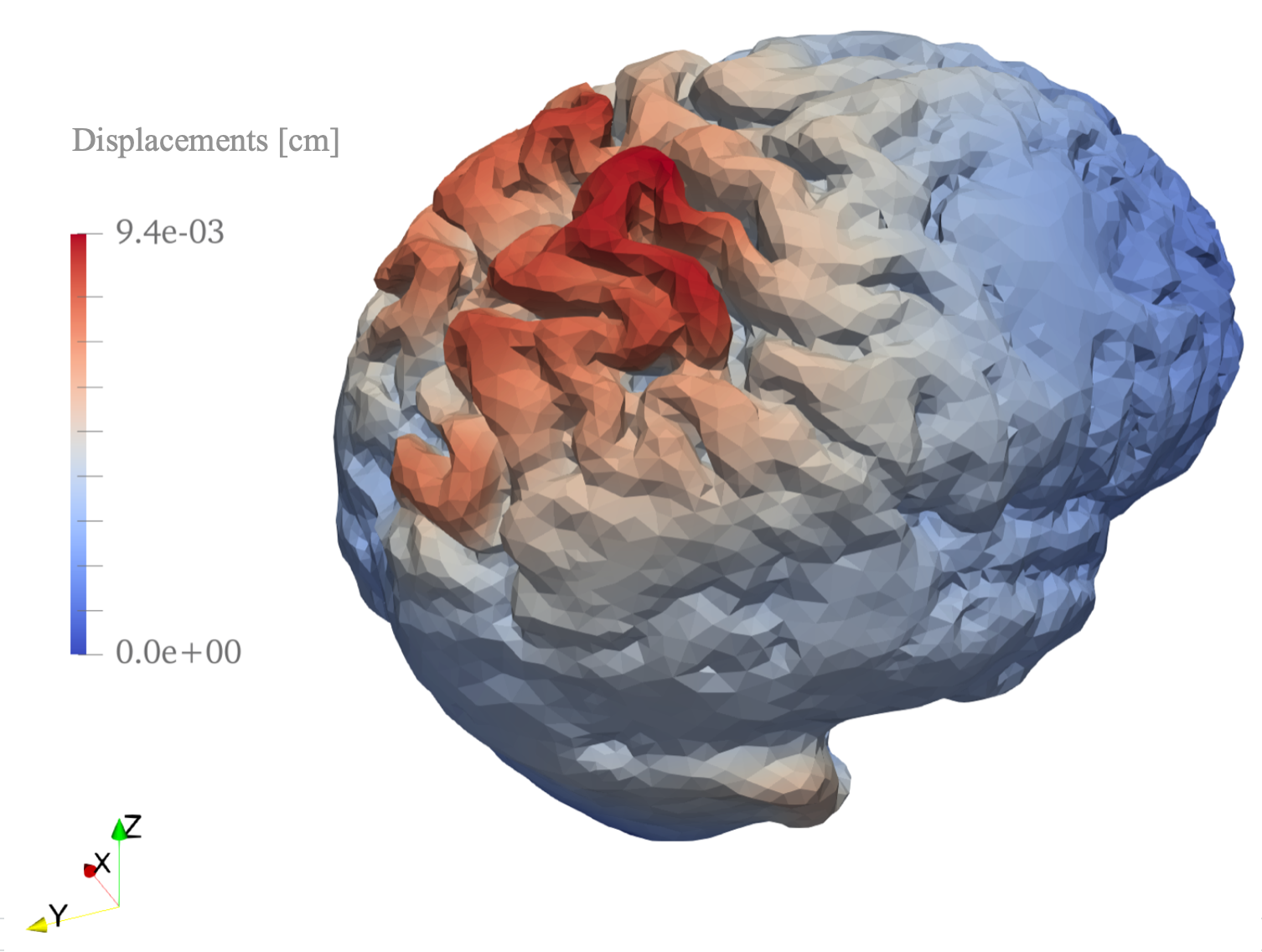}}
\subfigure[View of the displacement field on a horizontal cross section]
    {\includegraphics[height=5cm]{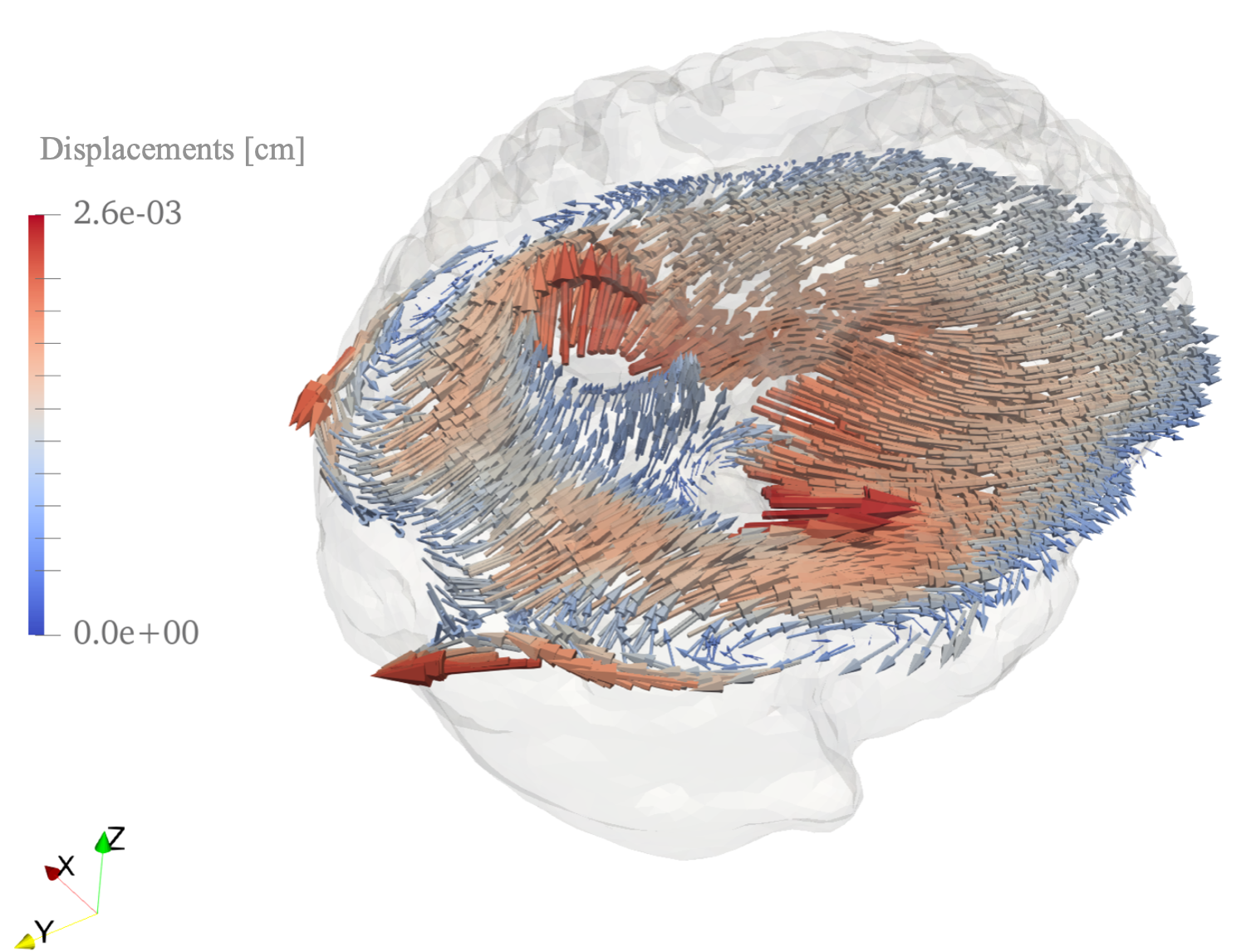}}
\subfigure[View of the displacement field on a vertical cross section]
    {\includegraphics[height=5cm]{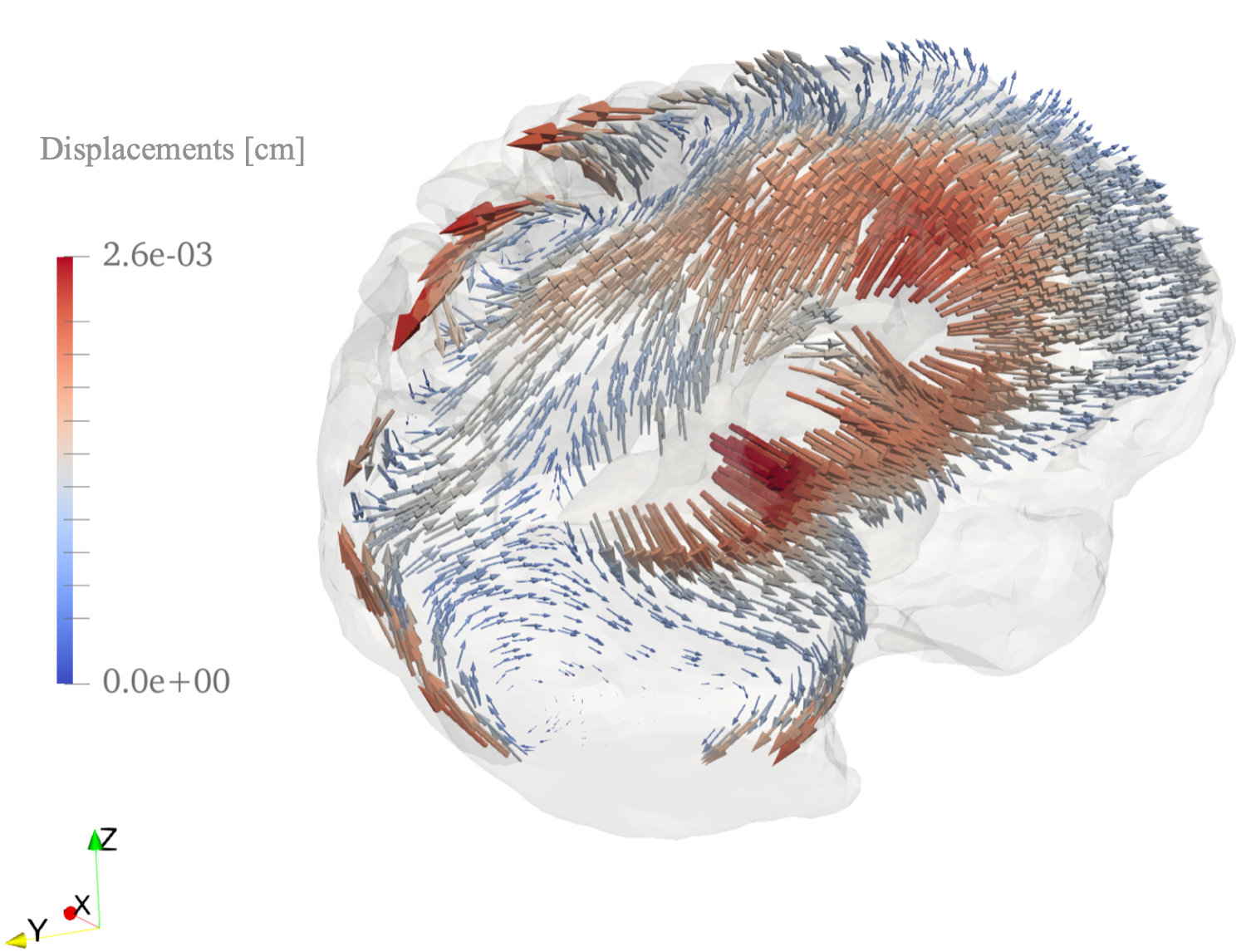}}
\subfigure[View of the pressure field]
    {\includegraphics[height=5cm]{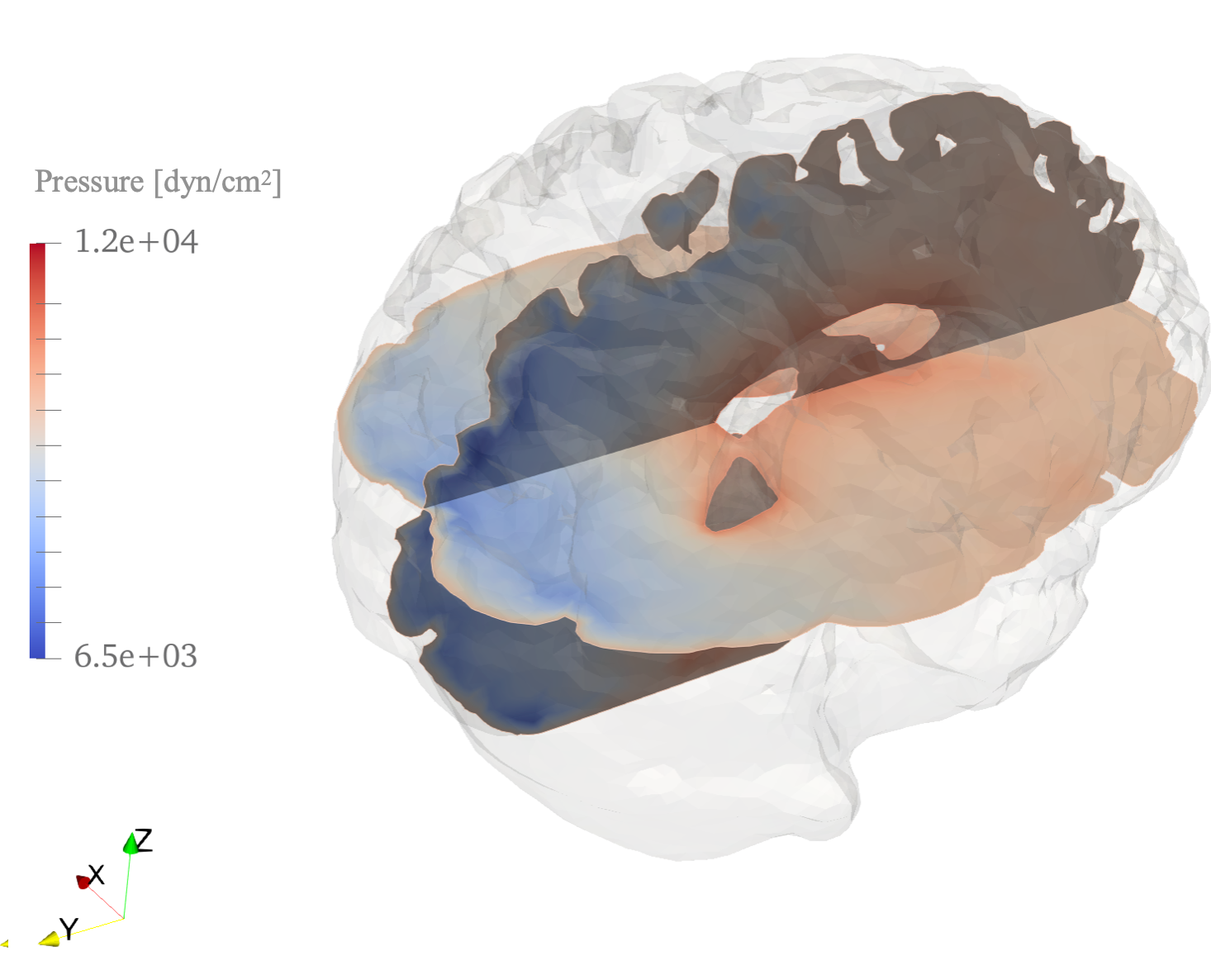}}
\caption{Snapshot of the numerical solution of \eqref{eq:the_model_PEDup_bc} for $\nu = 0.4$, $E=10^5$ dyn/cm\textsuperscript{2} and $\kappa = 10^{-8}$ cm\textsuperscript{2}, 15 ms after pulse is initiated.}
\label{fig:brain-forward}
\end{figure}

The training set manifold is defined by the outputs of 512 simulations, each including a time series of 40 time steps (per pulse cycle).
Thus, the training manifold is composed of a total of $K = 20480$ snapshots, each one being a finite element function of dimension $\cN \approx 131500$, four (degrees-of-freedom per nodal point) times the number of nodes of the FE mesh.
As an example, one of the snapshots in $\cM^{\text{training}}$ is shown in Figure \ref{fig:brain-forward}. 
The computations were carried out in parallel on machines with 768 GB of RAM and up to 72 threads.

\begin{figure}[!htbp]
\centering
\subfigure[View of the first four displacements modes]{\includegraphics[height=5cm]{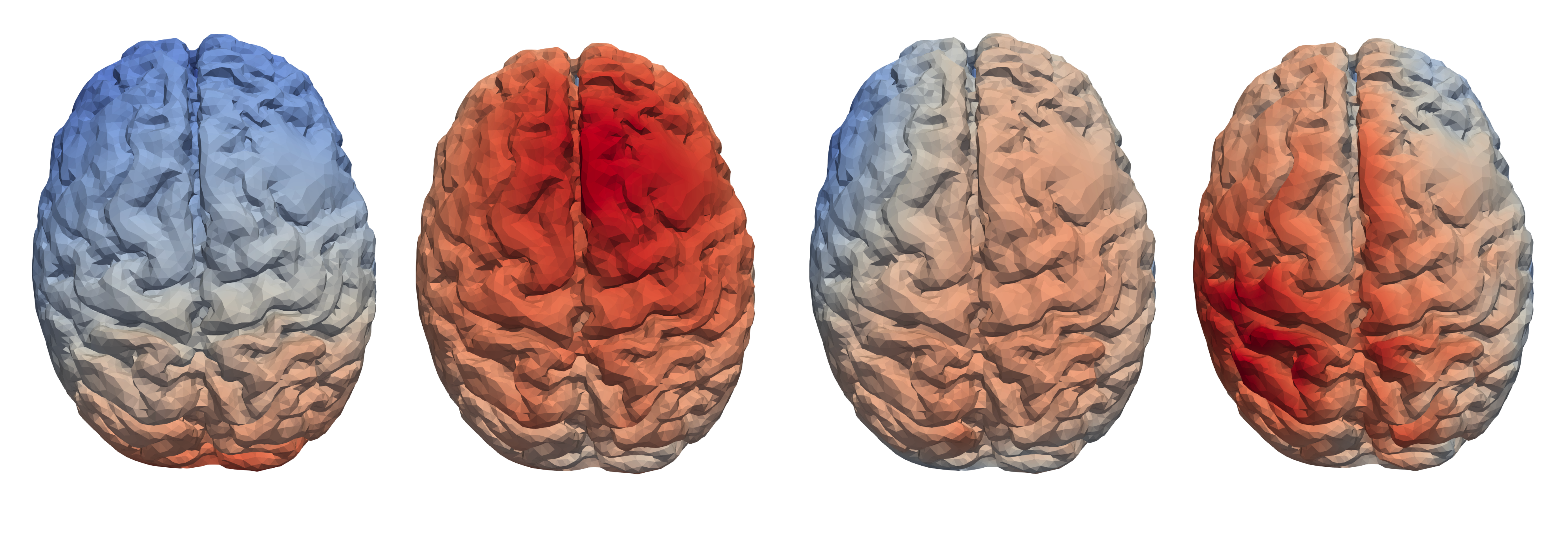}}
\subfigure[View of the first four pressure modes]{\includegraphics[height=3.8cm]{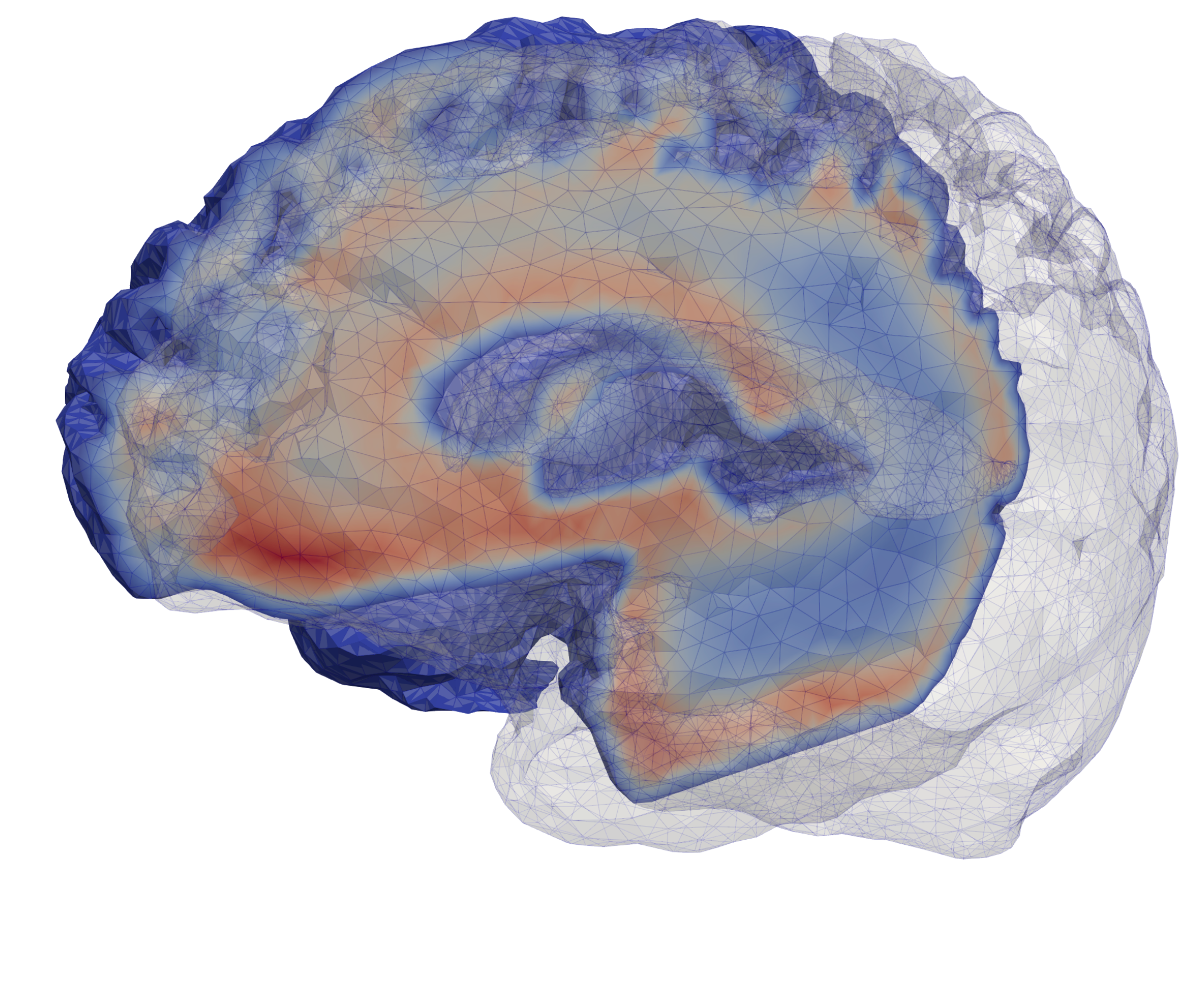}}
\caption{Basis functions (non-dimensional) of the reduced space $V_n$ for displacement and pressure. The functions are normalized in the $L^2$-norm.}
\label{fig:pod_modes_up}
\end{figure}

Subsequently, the snapshots are stored in a data matrix: $A \in \bR^{\cN \times K}$.
The reduced basis defining the space $V_n$ is then computed using the PCA of the snapshot matrix, as described in Section \ref{ssec:training}.
The first four modes of both displacement and pressure fields are depicted in Figure \ref{fig:pod_modes_up}.

\subsection{Joint reconstruction of displacements and pressure}

\review{As underlined above, 
proposed data assimilation approach aims to characterize 
pathological pressure gradients through MRE. To be able to infer pressure-dependent quantities, it is necessary to apply the PBDW method for the reconstruction of both the displacement and the pressure fields.}

The main challenge of this joint reconstruction strategy is a consequence of the fact that only displacement data are available, and hence the pressure state is thus \textit{invisible} to the data assimilation algorithm.

\rreview{The joint reconstruction problem can be formulated by considering}
a state space which is the direct product of two Hilbert spaces, i.e., $V_h = U_h \times Q_h$, with the induced scalar product
\begin{equation}
\< \( u, p \), \( v, q \) \>_{\zeta,V_h} \coloneqq 
\< u, v \>_{U_h} + \zeta \< p, q \>_{Q_h}
\,,
\label{eq:joint_hilbert}
\end{equation}
and the corresponding norm $\norm{(u,p)}_{\zeta} := \norm{u}_{U_h} + \zeta \norm{p}_{P_h}$.
The parameter $\zeta$ can be tuned in such a way to weight equally the norms of both components, and it can be
also seen as a conversion between the different physical units. 
This step enhances the stability of the computation of the reduced-order model. As it will be detailed below,
the scaling can also improve the quality of the stability constant $\beta(V_n, W)$. 
In practice, this parameter can be defined via
\begin{equation}\label{eq:zeta-up}
\zeta = \dfrac{\max_{u \in \manit}{\norm{u}_{V_h}}}
              {\max_{ p \in \manit}{\norm{p}_{V_h}}}.
\end{equation}

The fact that only displacements are observable can be formalized considering 
an observation space (see Section \ref{ssec:W-space}) of the form
\begin{equation}\label{eq:W-joint}
W= \text{Span}(w_1,\hdots,w_m) \times \{ 0 \} \rreview{:=W_U \times \{ 0 \}},
\end{equation}
where $w_i$, $i=1,\hdots,m$ are the Riesz representers of the functionals
$l_i: U_h \to \mathbb R$, $i=1,\hdots,m$,
defined on the displacements space $U_h$.

\rrreview{
\begin{remark}[Coupling of displacement and pressure modes]\label{rem:coupled}
In the case of the joint reconstruction, the reduced-order space $V_n$ 
is spanned by $n$ POD modes on the joint state, i.e., of the form
$$
\rho_i = (\rho_i^u,\rho_i^p),\, i=1,\hdots,n\,,
$$
where $\rho_i^u$ and $\rho_i^p$, $i=1,\hdots,n$, span
$n$-dimensional reduced-order models in $U_h$ and $Q_h$, respectively.

Notice that this construction enforces a coupling between the displacement and pressure components on the
reduced-order space $V_n$.
In particular, $V_n$ is not a tensor product of two reduced-order spaces for displacement and pressure.
\end{remark}

\begin{remark}[On the observability of the joint reduced-order space]
Assume that the condition \eqref{eq:pbdw_bound} is satisfied for a displacement 
reconstruction problem on $U_n = \text{Span}\left\{\rho_i^u,i=1,\hdots,n\right\}$,  with the given (displacement) observation space $W_U$, i.e., 
$$
\inf_{u \in U_{n} } \dfrac{ \Pi_{W_U}(u)}{\norm{u}_{U_h}}>0.
$$
From the coupling of the modes (Remark \ref{rem:coupled}) if follows that the linear operator
$$
\begin{aligned} 
& V_n & \to & \;U_n \\
v = & \sum_i \alpha_i (\rho_i^u,\rho_i^p) &  \mapsto & \;u_v : = \sum_i \alpha_i \rho_i^u
\end{aligned}
$$ 
is invertible, thus that there exists a 
constant $C>0$, depending on the basis elements, such that $\norm{v}_{\zeta} \leq C \norm{u_v}_{U_h}$.

Hence, 
\begin{equation}\label{eq:pbdw-joint-bound}
\beta\left(V_n ,W_U \times \{0\}; \zeta \right) = \inf_{v = (u_v,p_v) \in V_n}
\frac{ \Pi_{W_U}(u_v)}{\norm{v}_{\zeta}^2} >0\,,
\end{equation}
%
%
%
%
i.e., the existence of the solution
in the joint setting is guaranteed when a solution in the displacement setting exists.
In other words, the fact that $U_n$ does not contains elements 
which are orthogonal to $W_U$ implies, due to the coupling of the modes in the joint
reduced-order space,  the observability of the joint space $V_n$. 
\end{remark}

\begin{remark}[On the joint stability constant]
It holds (see also \cite{GLM2021}):
\begin{equation}\label{eq:pbdw-joint-bound-1}
\beta\left(V_n ,W_U \times \{0\}; \zeta \right) = \inf_{v = (u_v,p_v) \in V_n}
\frac{ \Pi_{W_U}(u_v)}{\norm{u_v}_{U_h} + \zeta \norm{p_v}_{P_h}}
\leq \underbrace{\inf_{u \in U_{n} } \frac{ \Pi_{W_U}(u)}{\norm{u}_{U_h}}}_{\beta \left( U_{n},W_U\right)}\,.
\end{equation}
Hence, the stability constant for the joint state estimation problem is lower than in the case of a pure displacement reconstruction problem.
\end{remark}
}

%

\begin{figure}[htbp]
\centering
\includegraphics[width=.3\textwidth]{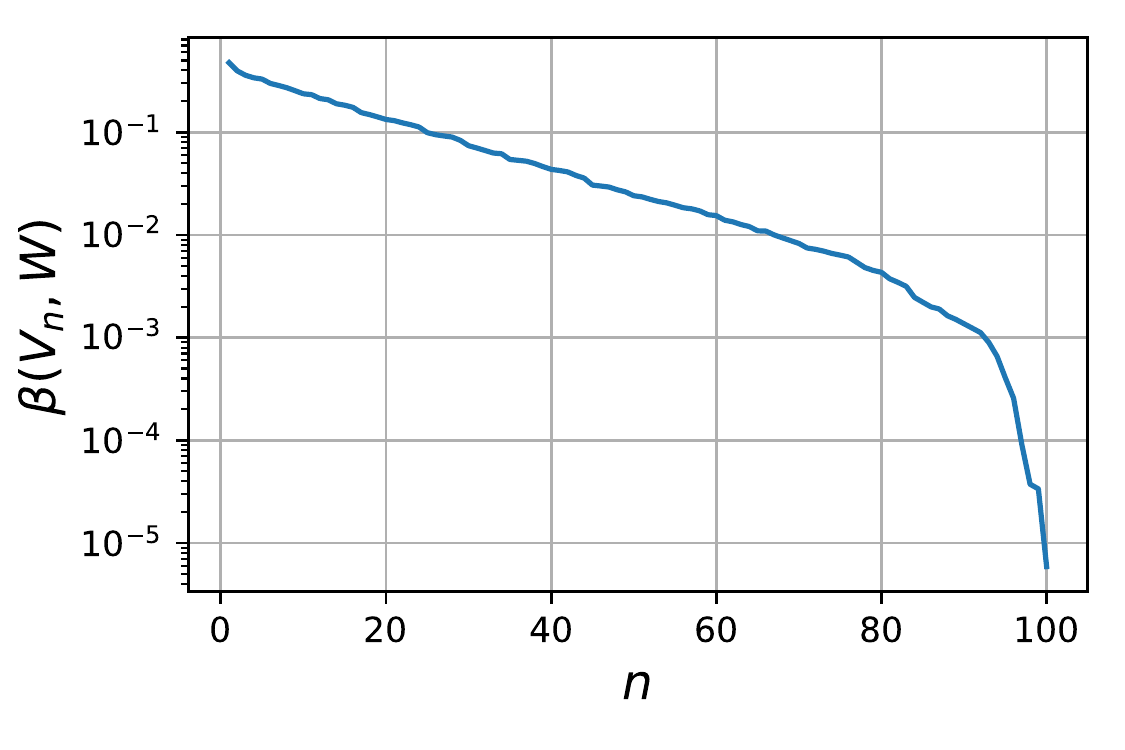}
\includegraphics[width=.3\textwidth]{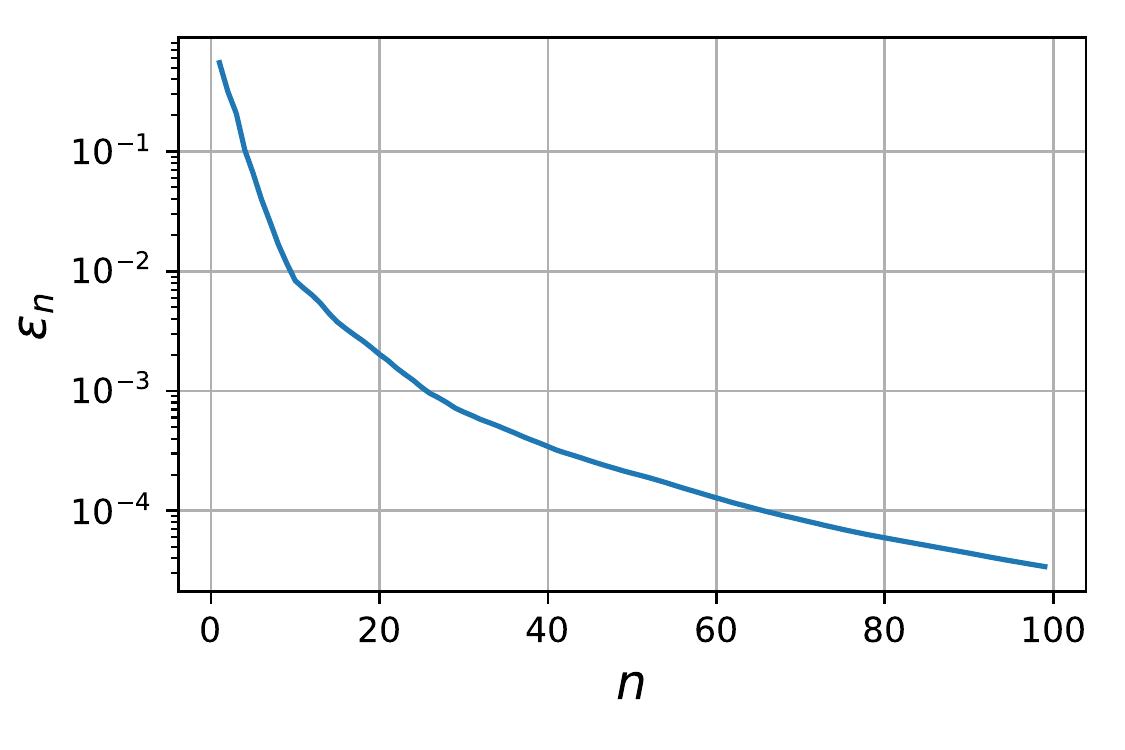}
\includegraphics[width=.3\textwidth]{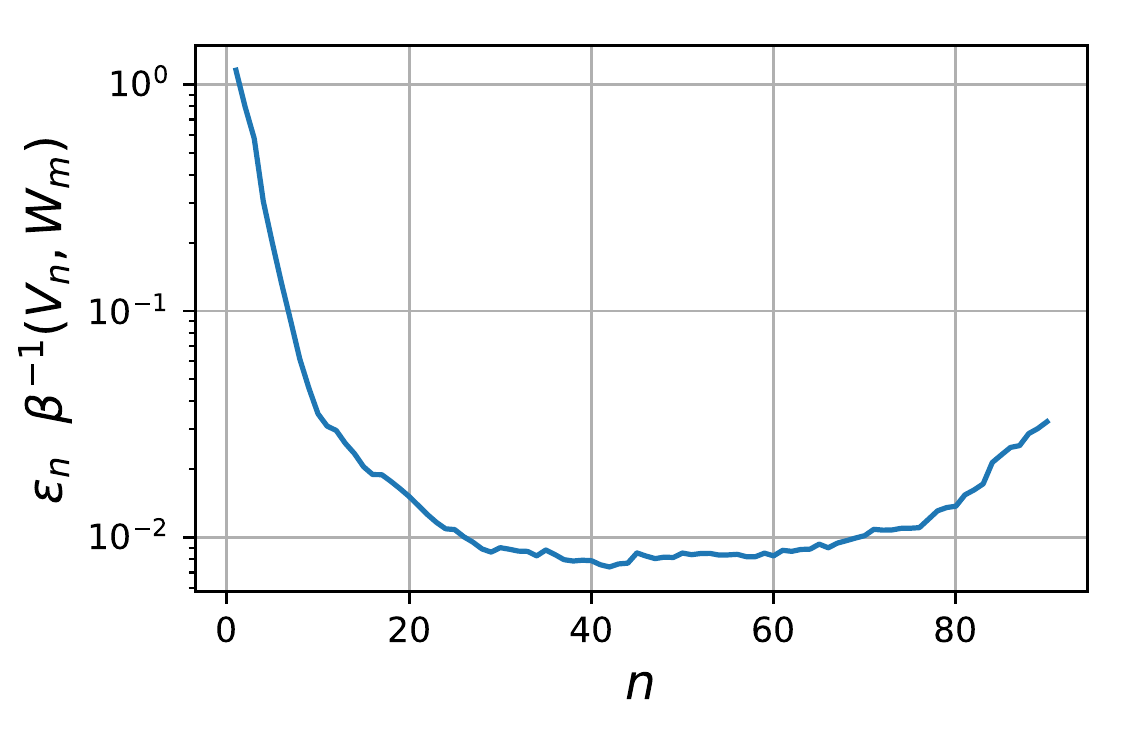}
\caption{
Left. Stability constant $\beta(V_n,W)$ \eqref{eq:pbdw_bound} (stability constant) between model and \review{measurements}, as a function of the reduced space dimension $n$. 
Center. Behavior of the \rreview{reduced-order model approximation} error $\epsilon_n$ \eqref{eq:mor_tails} as a function of the reduced space dimension. Right. Normalized \textit{a priori} error bound \eqref{eq:error-joint} as a function of the dimension of the reduced space dimension $n$.}
\label{fig:rec_apriori}
\end{figure}

As pointed out in Remark \ref{rem:beta}, the positivity of the stability constant can be assessed a priori.
For the considered case, the computed values of the stability constant as a function of the dimension $n$ of the 
\rrreview{joint reduced-order space $V_n$}
are depicted in Figure \ref{fig:rec_apriori}, left.

Using the joint state space means that, for $(u_{\text{true}},p_{\text{true}})  \in V_h$, the following error bound holds
\begin{equation}\label{eq:error-joint}
\begin{aligned}
\norm{(u_{\text{true}},p_{\text{true}}) - (u,p)}_{V_h}^2 & = \norm{u_{\text{true}} - u}^2 + \zeta \norm{p_{\text{true}} - p}^2
\\
& \le \left(\hat{\epsilon}_{U_h,n}^2 + \zeta \hat{\epsilon}_{P_h,n}^2\right) \beta^{-2} 
\left( \rrreview{V_n}, 
W_U \times \{0\}; \zeta \right)\,.
\end{aligned}
\end{equation}
\rreview{The model error $\hat{\epsilon}$ can be computed from the singular values of the snapshot matrix
(Figure \ref{fig:rec_apriori}, left), while the curve for the resulting error bound \eqref{eq:pbdw-error-bound} as a function of the dimension of the reduced-order space is shown in Figure \ref{fig:rec_apriori} (right).}
Based on this curve, the dimension of the reduced space $n$ can be chosen as the value that minimizes the quantity 
$\hat{\epsilon}_n\,\beta^{-1}$ ($n=41$ in the considered case).

\begin{remark}
The rescaling of the norms through the factor $\zeta$ \rreview{does not affect the well-posedness of the problem, but it} can improve the joint stability constant and, therefore, 
the error bound for the joint reconstruction. At the same time, \eqref{eq:error-joint} shows that this might yield to worse control on the
different components.
\end{remark}

\section{Results}\label{sec:results}

\subsection{Validation}\label{sssec:validation}
\label{sec:res_1}
\review{
The proposed MRE data assimilation algorithm is validated by reconstructing the poroelastic solution in $N_{\text{test}}=18$ test cases 
not included in the 512 samples used for generating the training manifold.
The 18 cases have been generated sampling the values of the physical parameters $\(\kappa, E, \nu, p_{\text{ventricles}}\)$
from a uniform distribution within the range of equation \eqref{eq:theta-training}, but discarding those
which were included in the original training set.
}
For each test case, the displacement field on the selected plane has been observed from the forward solution to generate the synthetic measurements.
Next, the overall displacement and pressure fields have been reconstructed using the reduced trained manifold space.
\review{
Finally, the reconstruction errors are evaluated for the joint state
\begin{equation} \label{eq:errors_L2_up}
e_{up,\theta_i}(t) = \frac{\norm{\(u_{\theta_i, \text{true}}, p_{\theta_i, \text{true}} \) - \(u_{\theta_i}^*, p_{\theta_i}^*\)}_{\zeta,L^2(\Omega)}}{\norm{
\(u_{\theta_i, \text{true}}, p_{\theta_i, \text{true}} \)}_{\zeta,L^2(\Omega)}},
\end{equation}
as well as for the single components, i.e., 
}
\begin{equation}
e_{u,\theta_i}(t) =
\frac{\norm{u_{\theta_i,\text{true}} - u_{\theta_i}^*}_{\review{L^2(\Omega)}}}
{\norm{u_{\theta_i,\text{true}}}_{\review{L^2(\Omega) }}}
\,,
~~
e_{p, \theta_i }(t) =
\frac{\norm{p_{ \theta_i,\text{true}} - p_{\theta_i}^*}_{ \review{L^2(\Omega)}} }{\norm{p_{\theta_i,\text{true}}}_{ \review{L^2(\Omega)}} }
\,, 
\label{eq:errors_L2}
\end{equation}
where $i=1,\hdots,18$, $u_{i,\text{true}}$ and $p_{i,\text{true}}$ stand for the numerical solution of the forward problem for parameters $\theta_i$ (ground truth), and $u_{\theta_i}^*$ and $p_{\theta_i}^*$ are the corresponding PBDW reconstructions. In addition, we consider the means over the validation set:
\begin{equation}
e_{up} (t) = \frac{1}{N_\text{test}} \sum_{i=1}^{N_\text{test}} e_{up,\theta_i} (t), \quad  e_{u} (t) = \frac{1}{N_\text{test}} \sum_{i=1}^{N_\text{test}} e_{u,\theta_i} (t), \quad e_{p} (t) = \frac{1}{N_\text{test}} \sum_{i=1}^{N_\text{test}} e_{p,\theta_i} (t).
\label{eq:means_error}
\end{equation}

\review{
Besides the time dependent errors, we also introduce the time-averaged quantities, e.g., 
\begin{equation}\label{eq:error-time-up}
e_{up}^T := \int_{[0,T]} e_{up} (t) ~ \dt.
\end{equation}

In addition, validation tests will consider the reconstruction in the case of noisy measurement. To model noisy data,
we assume that the observations can be seen as linear functionals of the form
\begin{equation}\label{eq:l-noise}
\hat{\ell_i} \( u \) = \ell_i + \epsilon_{gn}^i,
\end{equation}
where $\epsilon_{gn}^i \approx \cN\(0, \sigma \)$. The standard deviation is chosen as
$$
\sigma = \Xi \max_{i,t} \ell_i(u(t,\cdot)),
$$
where $\Xi$ denotes the noise intensity relative to the maximal signal in time and space.
}

\begin{figure}[!htbp]
\centering
\includegraphics[height=6.5cm]{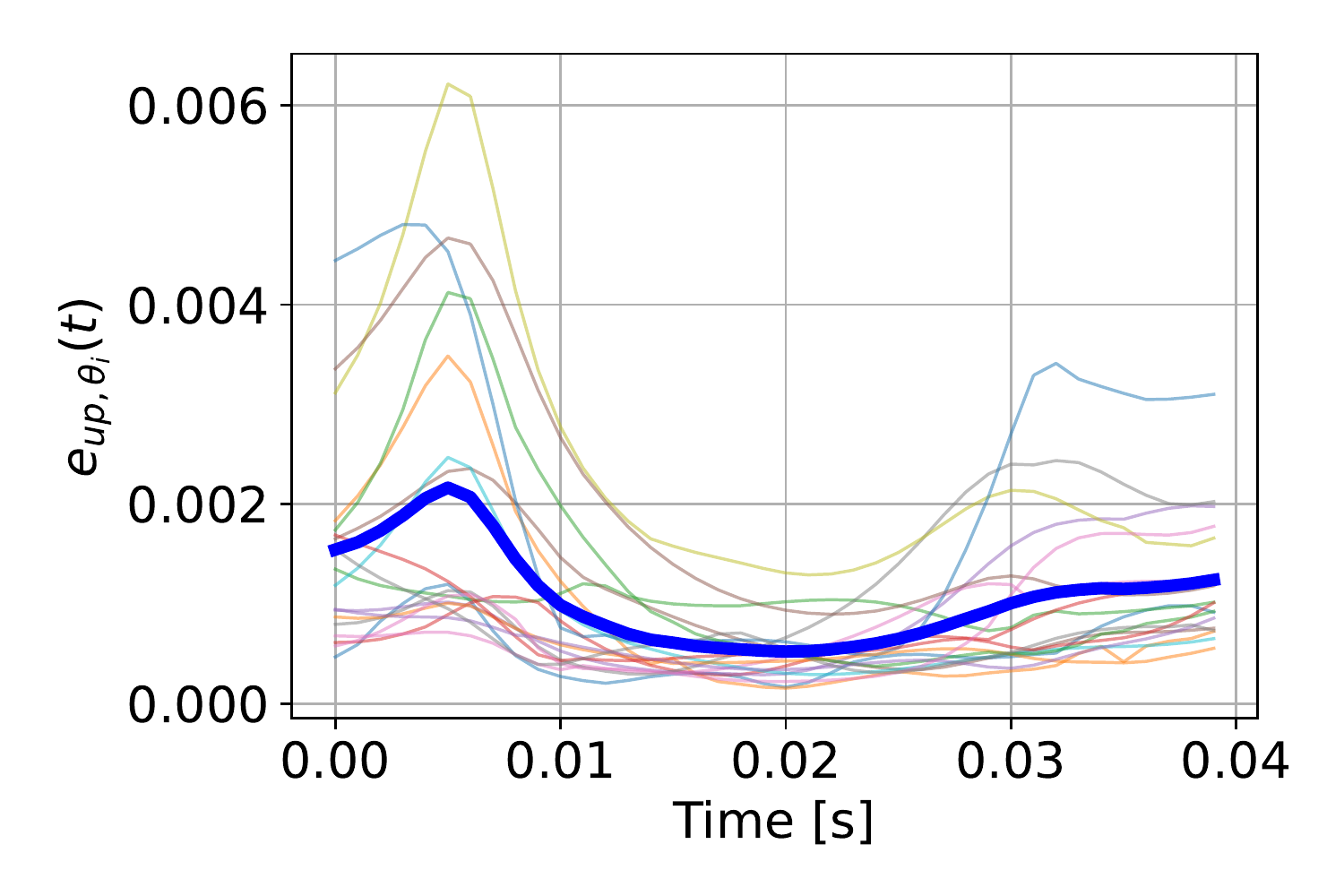}
\caption{Joint reconstruction error \eqref{eq:errors_L2_up} for the 18 different test cases used for validation. The blue curve shows the average error
over the 18 samples $e_{up} (t)$ from \ref{eq:means_error}.}
\label{fig:rec_joint}
\end{figure}

\begin{figure}[!htbp]
\centering
\subfigure[Displacements]{\includegraphics[height=5cm]{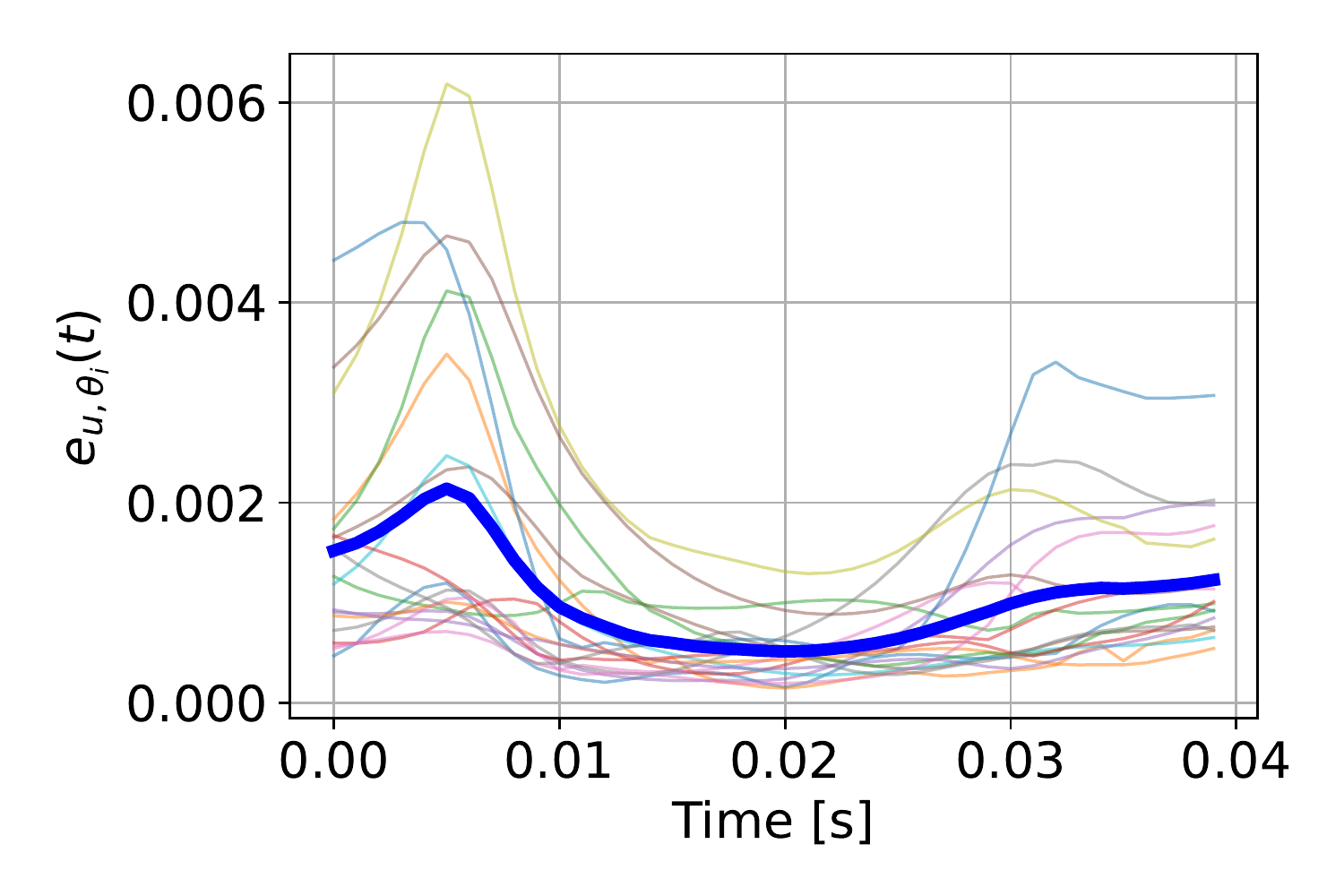}}
\subfigure[Pressure]{\includegraphics[height=5cm]{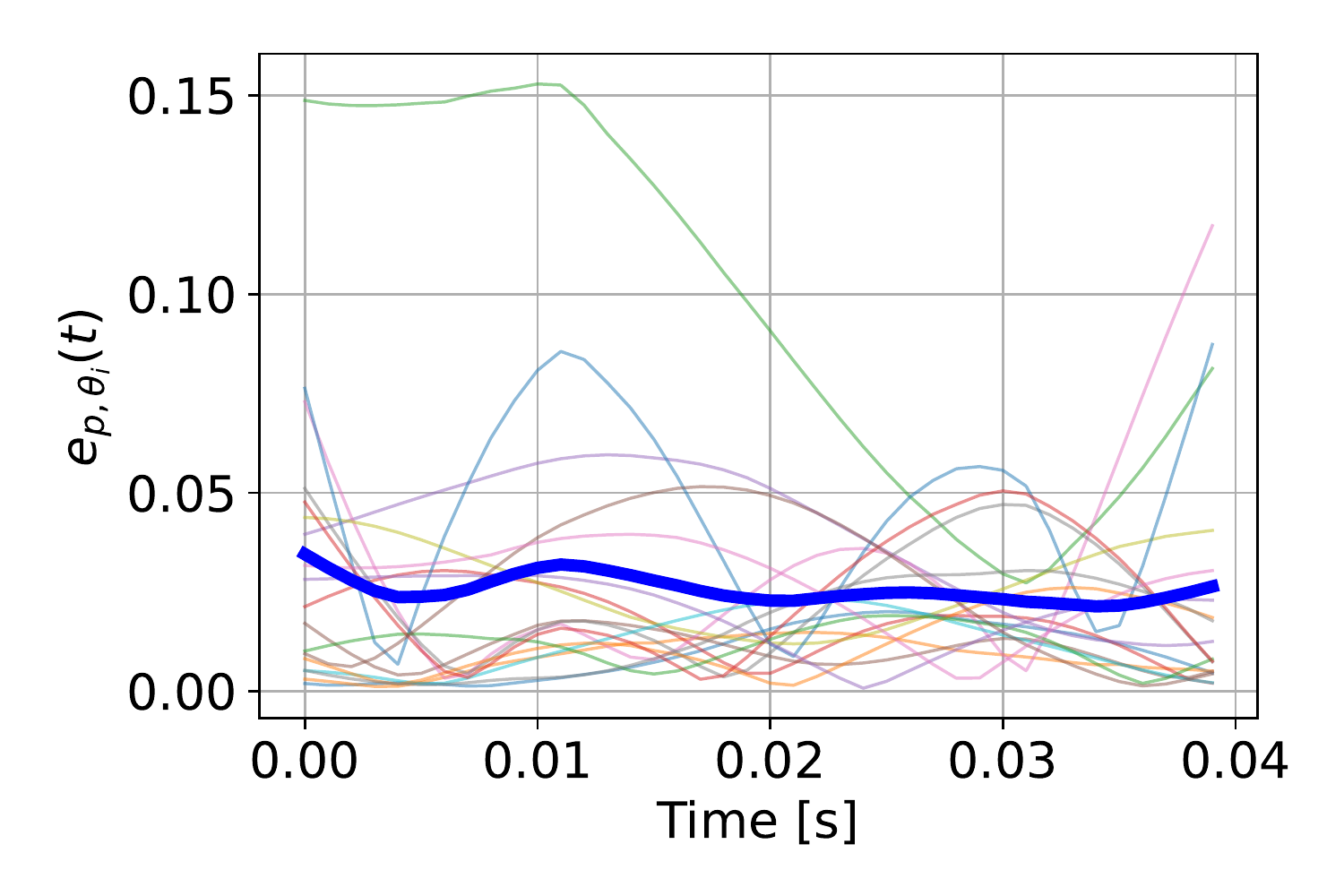}}
\caption{$L^2$ reconstruction errors from \eqref{eq:errors_L2} for displacement (left) and pressure (right) for the 18 samples used for validation. The blue curve shows the average error over the 18 cases from \eqref{eq:means_error}.}
\label{fig:rec_error}
\end{figure}

The joint reconstruction errors for all samples, using the metric introduced by the inner product \eqref{eq:joint_hilbert}, are shown in Figure \ref{fig:rec_joint}.
The reconstruction errors for each field considered (displacements and pressure) are shown separately in Figure \ref{fig:rec_error} (left and right, respectively).
In particular, it shows that the displacements are well reconstructed in all cases.

Moreover, we also observe that the pressure is reconstructed with satisfactory accuracy (peak errors mostly below 10\% over the whole cycle), with one single case showing a peak error of the order of 15\% (average error below the 5\%).

\begin{figure}[!htbp]
\centering
\subfigure[Average $e_{up}^T$]{\includegraphics[height=3cm]{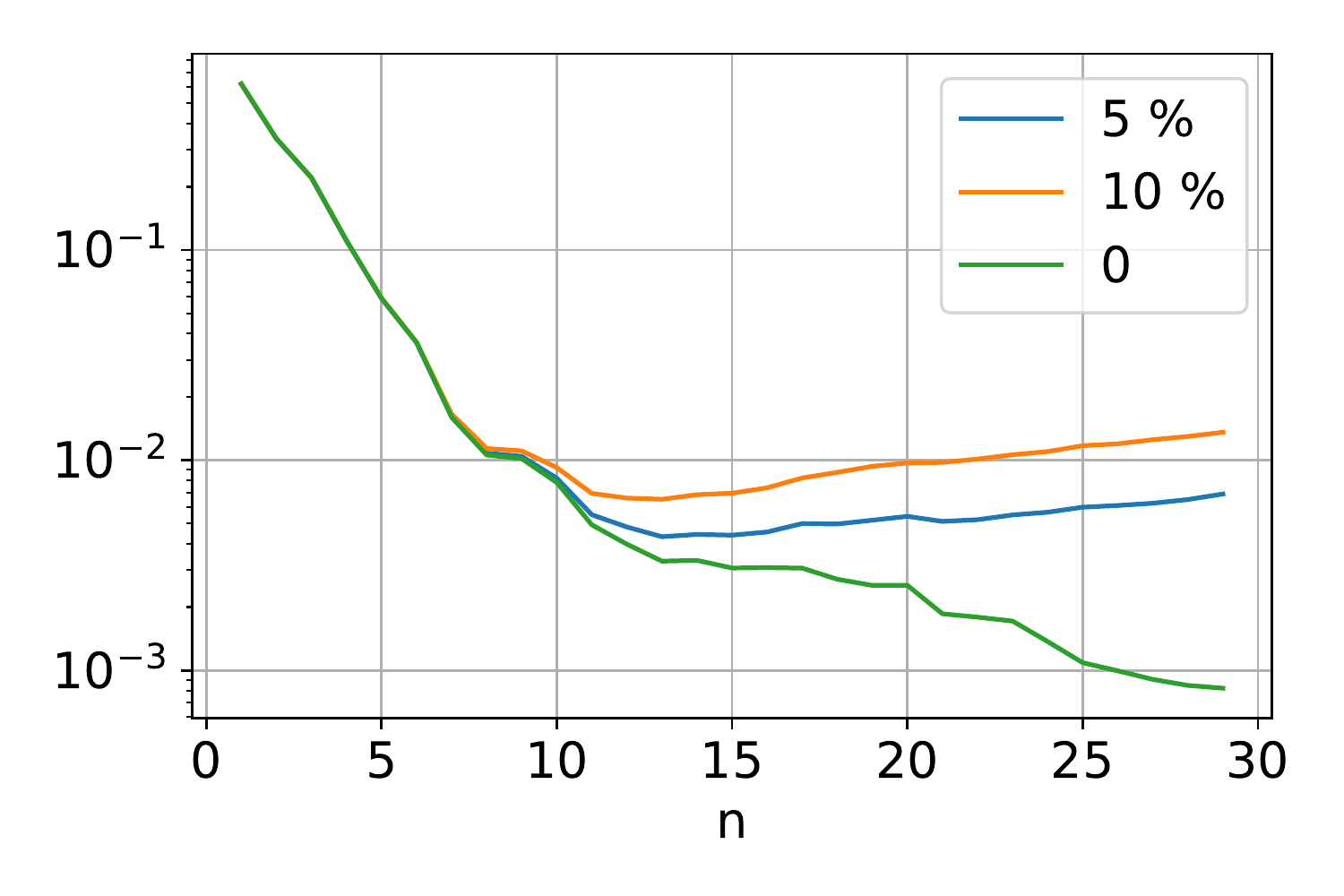}}
\subfigure[Average $e_{u}^T$]{\includegraphics[height=3cm]{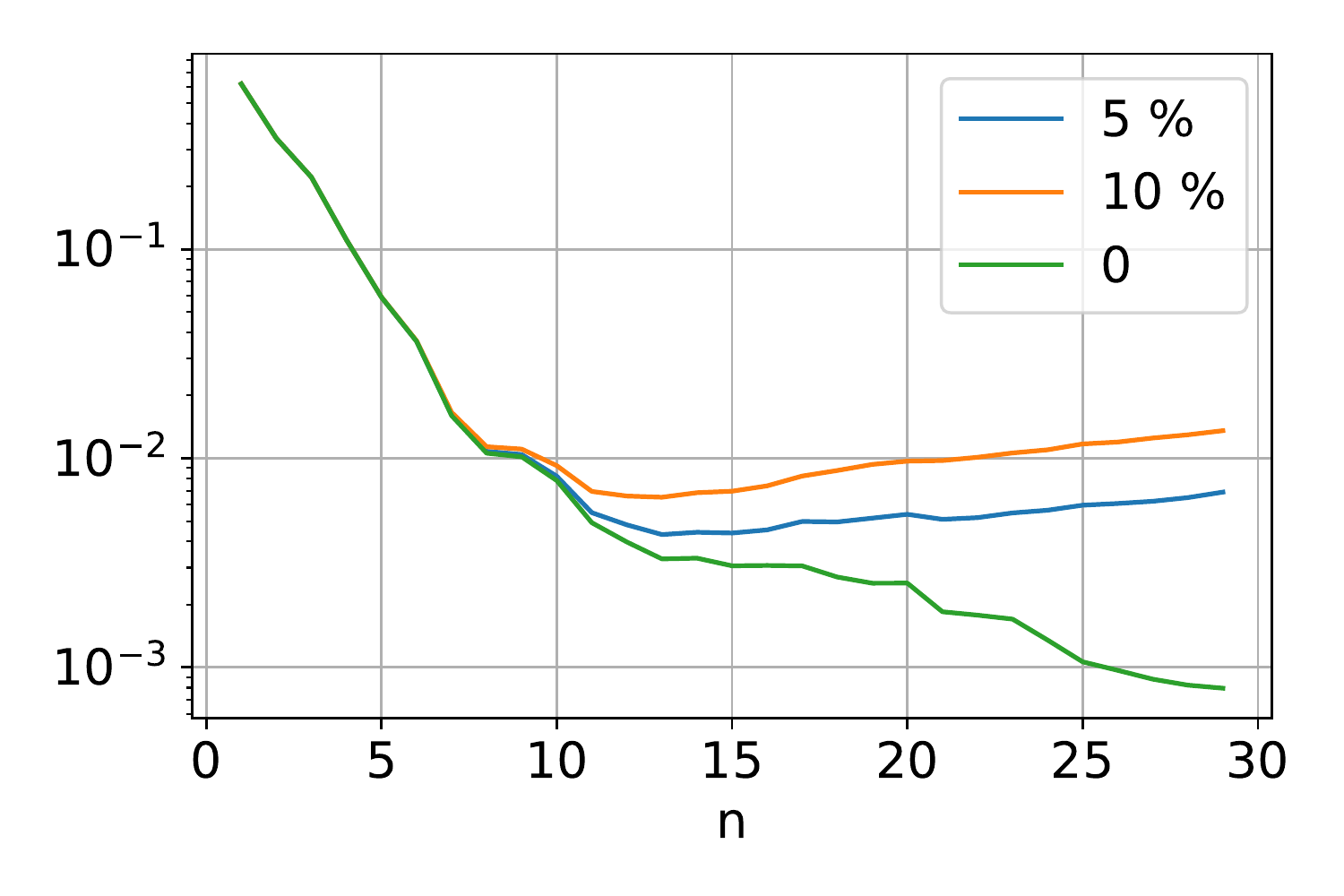}}
\subfigure[Average $e_{p}^T$]{\includegraphics[height=3cm]{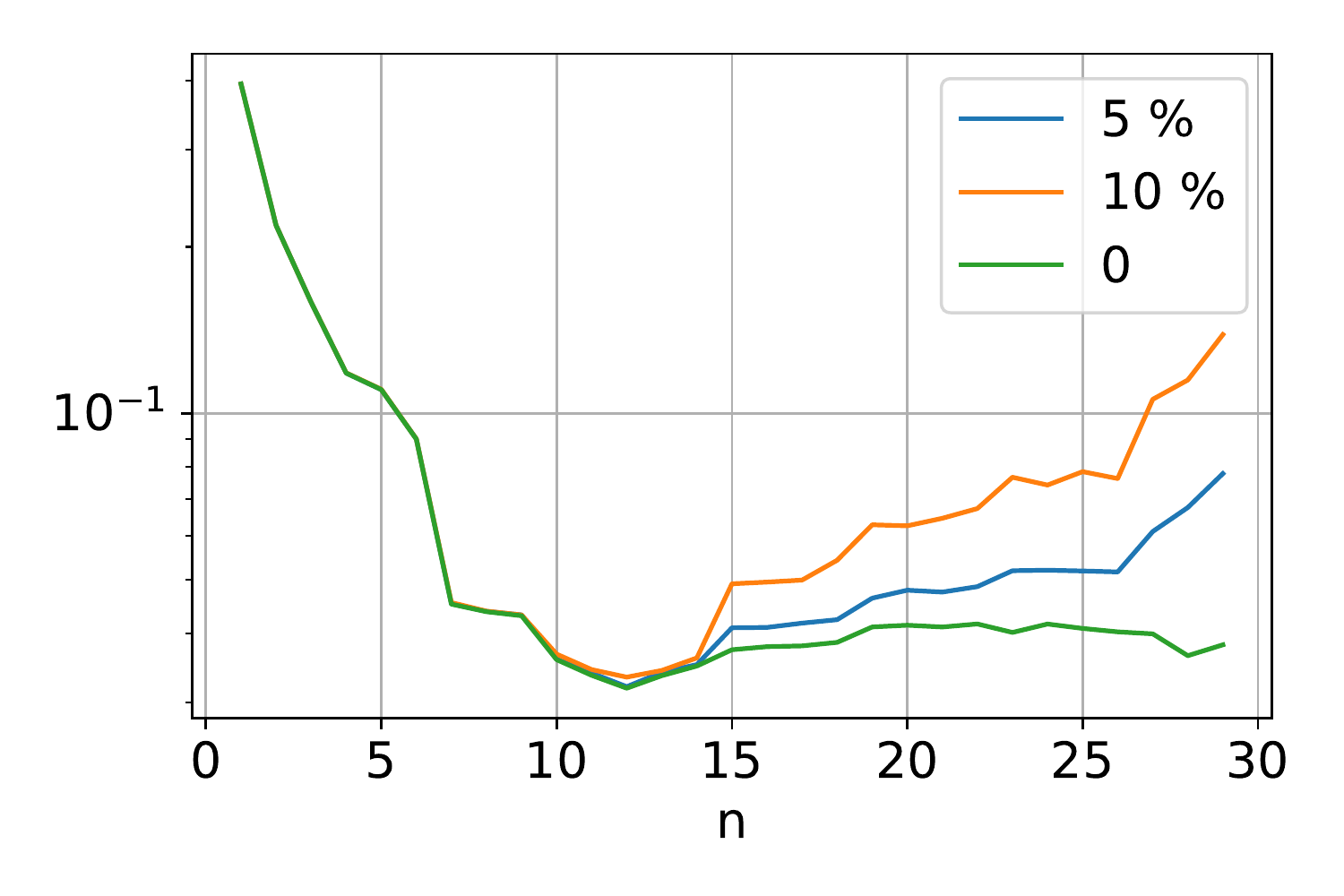}}
\caption{\review{Time-averaged reconstruction error for different levels of noise in the observations ($\Xi=\frac{1}{10}, \frac{1}{20}, 0$)as a function of the dimension of the reduced-order model $V_n$
(compare with Figure \ref{fig:rec_apriori}). Left: error for the joint-reconstruction. Center: displacement
error. Right: pressure error. }}
\label{fig:noise_vs_n}
\end{figure}

\review{A summary of results for the reconstruction error when using noisy data is shown in Figure \ref{fig:noise_vs_n}. 
The plot displays the values of the time-averaged error as a function of the reduced-order space dimension $n$.
One can observe that, the optimal value of $n$ is below the theoretical value obtained in the case of noise-free data ($n=41$).
This discrepancy arises from the fact that the bound \eqref{eq:pbdw_bound} is altered by the presence of noise. 
No theoretical results for the noisy case are available to the author's knowledge.
For a value of $n$ between 10 and 20, the joint reconstruction error stays below a reasonable level. On the other hand,
the error in pressure is considerably larger than the error in displacement, and it does not improve increasing
the dimension of $V_n$.

Moreover, one also sees that the pressure error, for the optimal value of $n$, is not considerably affected by the
presence of noise in the displacement data.
}

\begin{figure}[t]
\centering
\subfigure[Displacement field]{
\includegraphics[height=5cm]{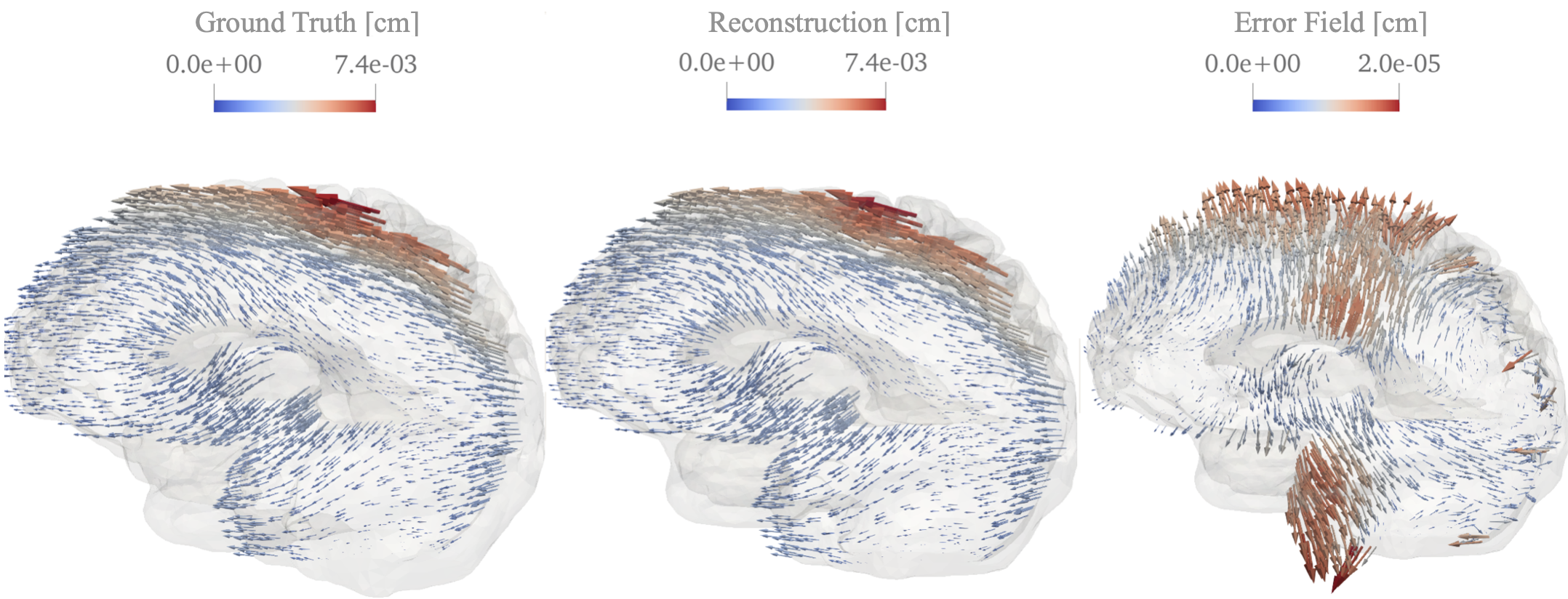}}
\subfigure[Pressure field]{
\includegraphics[height=5cm]{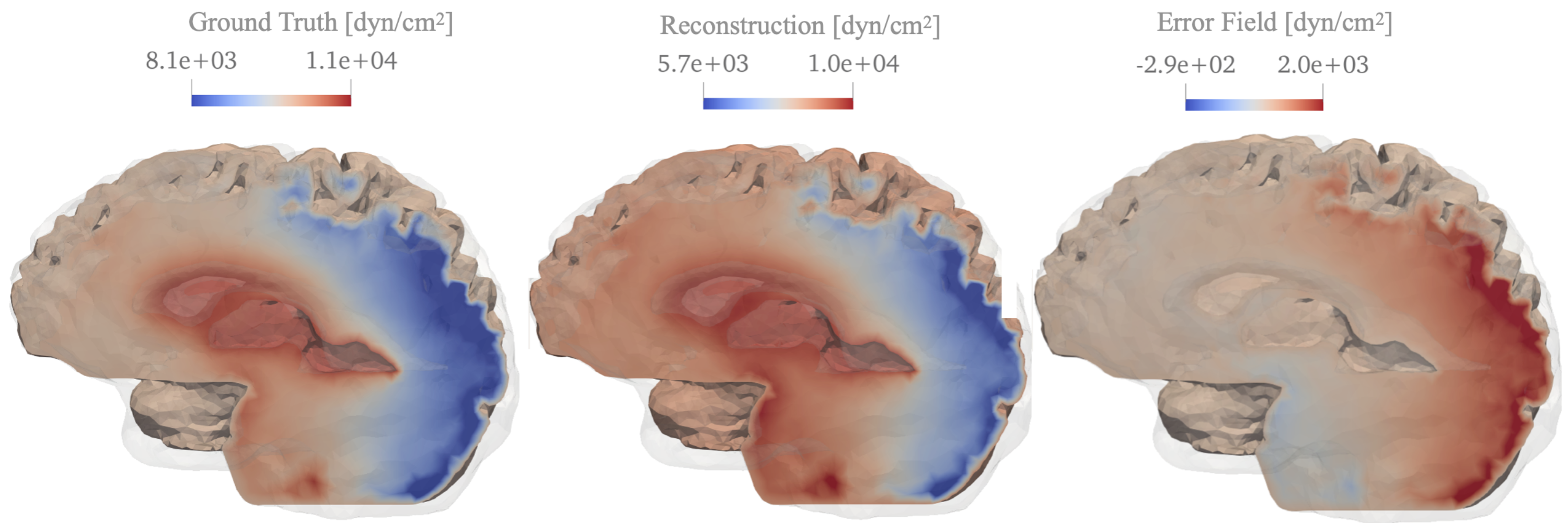}}
\caption{Reconstruction at $t=10$ ms after the pulse start, of displacement (top) and pressure (bottom) fields for the scenario with the largest reconstruction error among the 18 samples used for the validation \review{(noise-free measuremens)}.} 
\label{fig:rec_wc}
\end{figure}

A snapshot of the results in the scenario with the largest values for the model error (Equation \eqref{eq:mor_tails}) is visualized in Figure \ref{fig:rec_wc}, which compares the reference solution, the PBDW reconstruction, and the corresponding relative difference.

\review{
\subsubsection{Effect of increasing the amount of available data}
It is relevant to explore the impact of increasing the amount of available data for the reconstruction. 
For this purpose, we can increase the number $m$ by adding further slices to the set of observed voxels.
Namely, instead of using one single slice at $z=0$, we consider (i) the case of three slices and (ii) the 
case when the data over the full brain are available (at a lower resolution than the finite element solution).
The results for these tests are illustrated in Figure \ref{fig:up_vs_m}.

\begin{figure}[!htbp]
\centering
\subfigure{
\includegraphics[height=3.5cm]{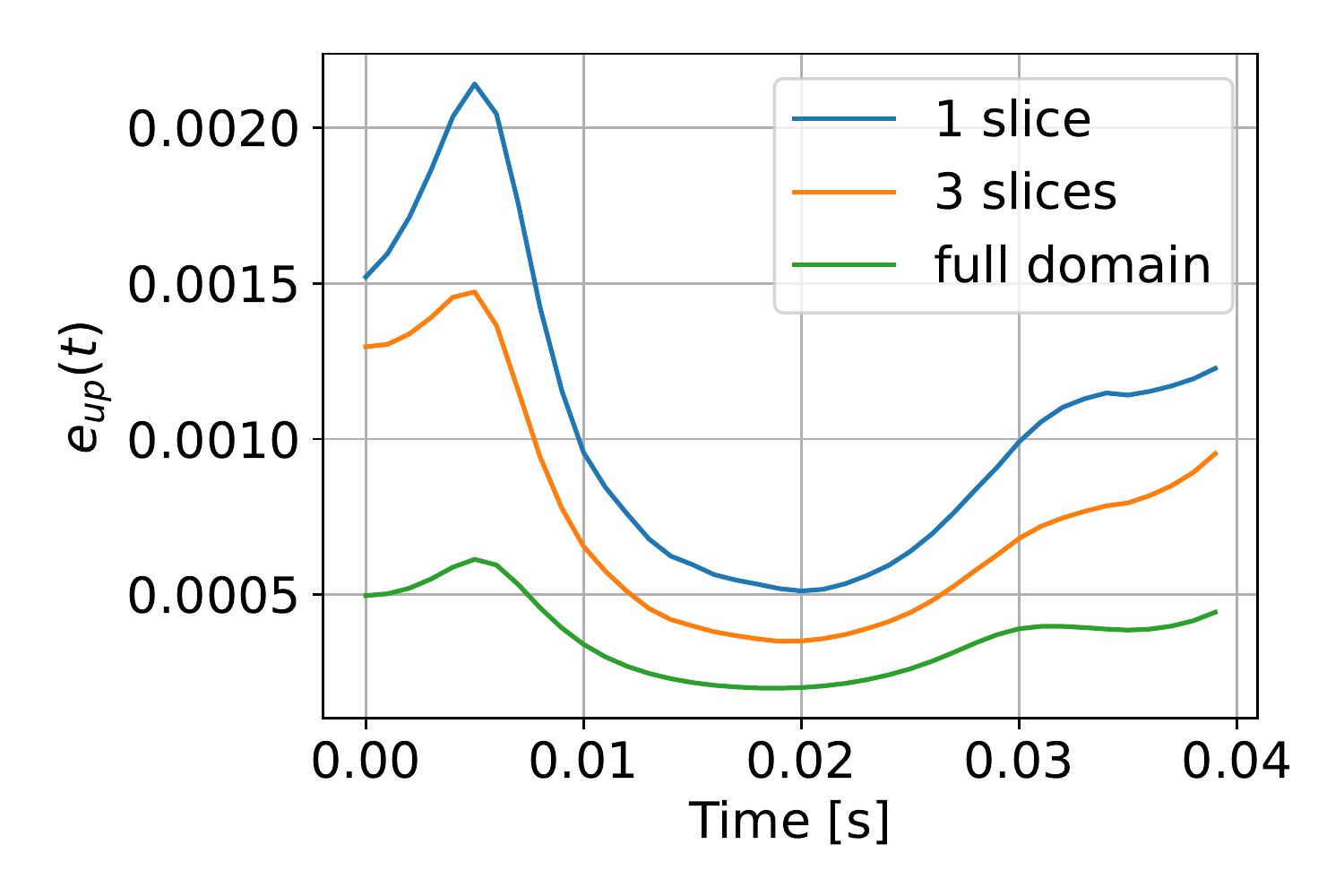}}
\subfigure{
\includegraphics[height=3.5cm]{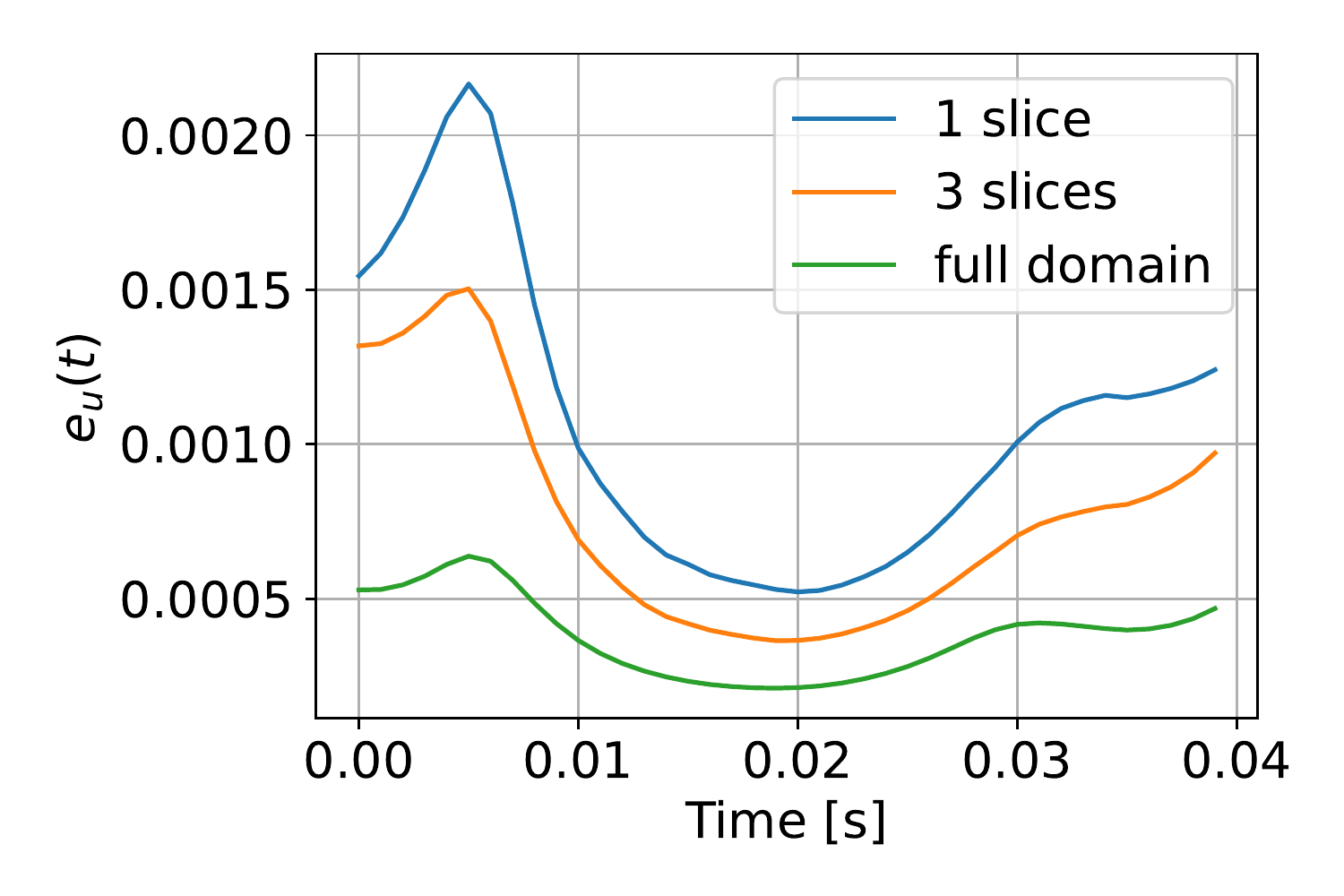}}
\subfigure{
\includegraphics[height=3.5cm]{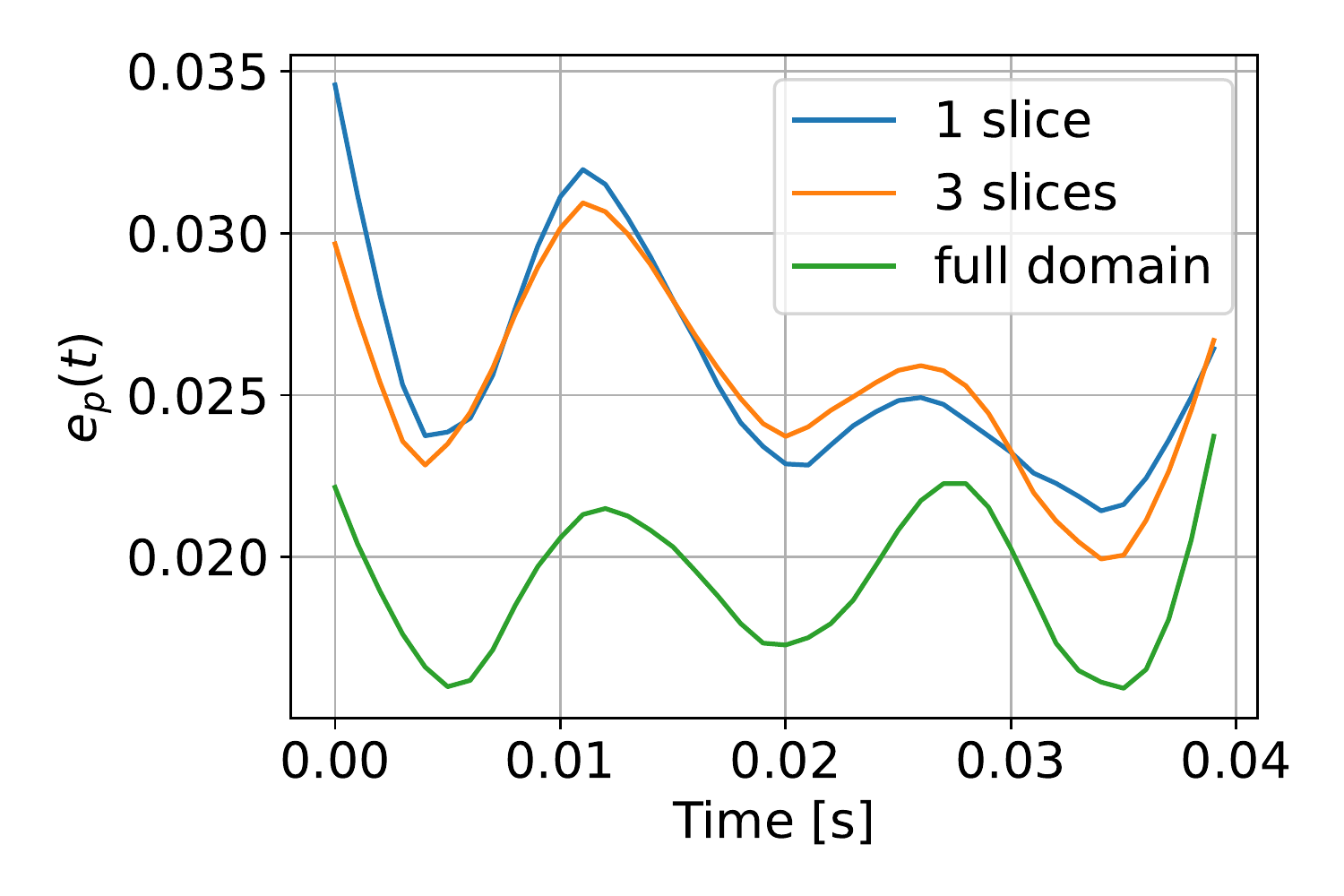}}
\caption{Reconstruction error for test cases with different data sets. The data set size is $m=600$ (200 voxels) for the single slice case, $m=2295$ (765 voxels) for the 3-slices case, and $m=5409$ (1803 voxels) for the full domain case.}
\label{fig:up_vs_m}
\end{figure}
One observes that the quality of the displacement reconstruction increases consistently with $m$, as it is expected to be theoretically. Uncoupling the errors (as defined at \eqref{eq:errors_L2}) one sees that the error in the pressure
decreases when the field is measured over the full brain, but the approximation 
is not improved when increasing the number of observed slices from 1 to 3.
These results hence demonstrate the relevance of a robust, physics-based, data assimilation algorithm in the case of non-observed
quantities.
}

\begin{figure}[!htbp]
\centering
\includegraphics[height=5cm]{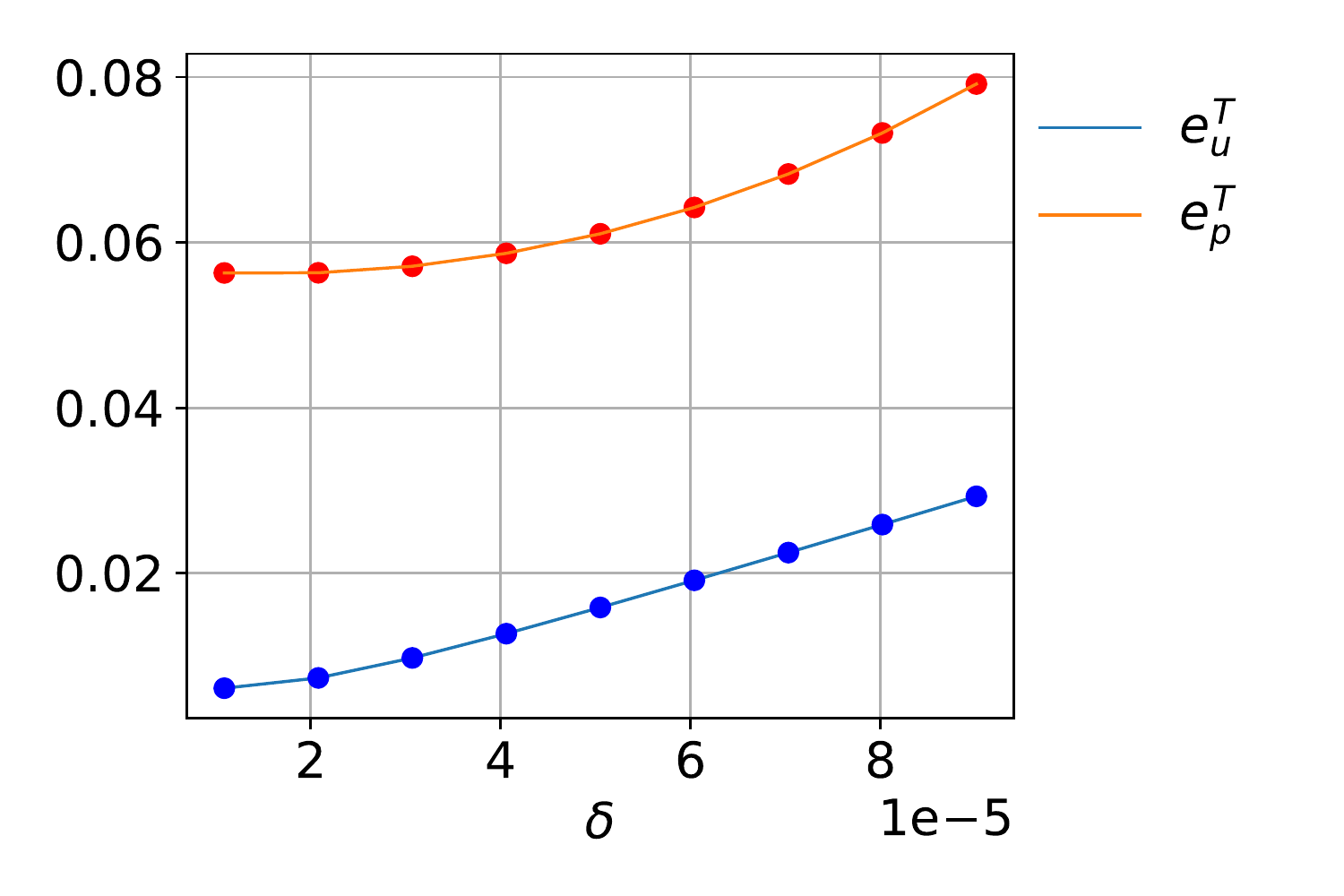}
\caption{\review{Reconstruction results taking into account model error (boundary conditions) on $\Gamma_{\text{neck}}$. 
The error stays reasonable until the magnitude of $\delta$ reaches the order of magnitude of the 
vibration on the boundary ($\approx 10^{-4}$ cms).}}
\label{fig:model_error_neck}
\end{figure}

\review{
\subsubsection{Model mismatch}
A further aspect, which is relevant for concrete applications, is the robustness of the reconstruction
in the case of a mismatch between the model used for the reconstruction and the physical model
from where the data are obtained.
From the point of view of the data assimilation framework, this mismatch corresponds to the case
in which the reconstruction belongs to an ambient space $\tilde V$, which is different from the space $V$ 
considered for building the solution manifold $\cM$ (and the corresponding reduced-order models $V_n$).
From the practical point of view, a model mismatch can be due to (i) a \textit{physical} mismatch, i.e., a mismatch 
between the PDE and the
underlying physical system, or (ii) a \textit{PDE-parameter} mismatch, i.e., an error in the parameters or boundary conditions defining the considered PDE model.
We focus on the latter case. In fact, it shall be noticed that mismatch at the level of the physical model (e.g., 
assuming a poro-viscoelastic mechanics, instead of a poroelastic only) should be handled by directly including
them in the manifold. Especially in the case of pressure reconstruction (a non observed quantity whose
approximation relies on the goodness of the approximation of the underlying physics), this type of mismatch would naturally
yield larger errors.

In what follows, we consider the presence of small errors on the boundary, which can model
a mismatch between the exact and patient geometry, as well as an error in the approximation of the exact 
Dirichlet boundary condition on the neck surface.
To this purpose, we generate the observations using a non-zero boundary condition on $\Gamma_{\text{neck}}$, of the form
$u_{\Gamma_{\text{neck}}} = \( \delta,\delta,\delta\)$, with a small $\delta > 0$.  The reconstruction is, however, sought in the same original manifold.
The results, shown in Figure \ref{fig:model_error_neck}, indicate that the approximation achieved by the data assimilation algorithm
is still satisfactory, as long as the boundary perturbation stays small, compared with the magnitude of the
overall displacement field.
}

\subsection{Non-invasive characterization of pressure-dependent biomarkers}\label{ssec:pv-compute}

In this paragraph, the characterization of elevated pressure gradients from displacement data of the data assimilation algorithm is interrogated.
For this purpose, the following quantity is introduced:

\begin{equation}\label{eq:pv}
p_{\text{v}}(p) \coloneqq 
\frac{1}{T} \int_{0}^T{\left[\frac{\int_{\Gamma_{\text{ventricles}}} p(x,t)~{\rm d}\Gamma}{\int_{\Gamma_{\text{ventricles}}}~{\rm d}\Gamma} \,\right] {\rm d}t}
\,,
\end{equation}
which represents the average, with respect to time and over the whole surface of the ventricular CSF pressure.
This indicator is important to evaluate the presence
of abnormal ICP in the brain, as well as characterize
increased ventricle volume.

The purpose of this numerical experiment is to test the accuracy of the ventricular pressure reconstruction (see Equation \eqref{eq:pv}), and interrogate the capability of the framework to stratify intracranial pressure (physiologic or pathological) by solely relying on limited displacement observations.
The numerical experiment has been designed as follows.
Two sets of synthetic individuals have been created, each set containing 8 individuals (16 in total). 
For each individual, the parameters $\kappa$, $E$, and $\nu$ have been sampled within the range specified in \eqref{eq:theta-training}.
The differentiation between physiological levels (normal) and pathological levels (increased) of ICP has been modeled by sampling the pressure boundary condition parameter within the aforementioned ranges.

\begin{equation}
p_{\text{ventricles}} \in \[ 1, 1.02 \] \times 10^4 \; \text{dyn/cm\textsuperscript{2}}
\,,
\label{eq:p-health}
\end{equation}
for the first group (normal intracranial pressure), and
\begin{equation}
p_{\text{ventricles}} \in \[ 1.08, 1.1 \] \times 10^4 \; \text{dyn/cm\textsuperscript{2}}
\,,
\label{eq:p-pathological}
\end{equation}
for the second group (increased intracranial pressure).

Hence, since the value of the outer CSF pressure is set to $10^4$ dyn/cm\textsuperscript{2}, the pathological cases (those with an increased pressure) are characterized by an ICP gradient between $800$ and $1000$ dyn/cm\textsuperscript{2} (0.8 to 1 cm-H\textsubscript{2}O), while the pressure difference for healthy scenarios (normal pressure) is up to 200 dyn/cm\textsuperscript{2} (0.2 cm-H\textsubscript{2}O), in line with the range of values obtained in the numerical experiments shown in \cite{li_etal_2014}.
For each of these 16 synthetic patients, the background PDE problem \eqref{eq:forward-weak}--\eqref{eq:A-weak} has been solved numerically, sampling the displacement data on a set of slices (see Figure \ref{fig:pbdw_measures}).
Then, the reduced-order model trained as discussed in Section \ref{sssec:validation}, is used to reconstruct the pressure field and compute the  ventricular pressure evaluating \eqref{eq:pv}.

\begin{figure}[!htbp]
\centering
\includegraphics[height=8cm]{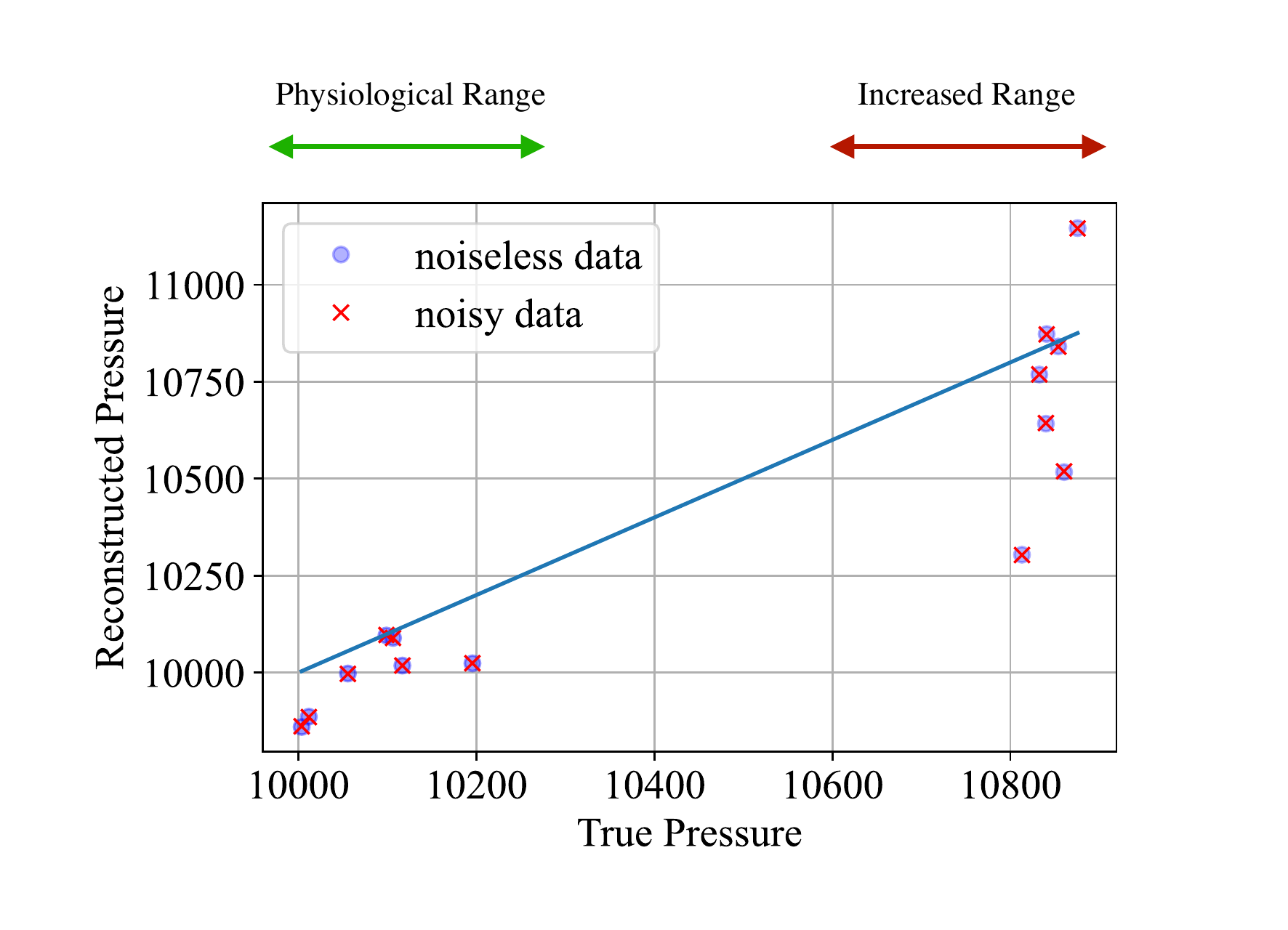}
\caption{Reconstructed ventricular pressure \eqref{eq:pv} versus the values computed using the reference values for the synthetic patients. Each scenario has been labeled as "normal" or "increased", depending on the value of $p_{\text{ventricle}}$.
The green dots correspond to the scenarios originally sampled from the normal range \eqref{eq:p-health}, whilst the orange dots correspond to the simulations with elevated pressures in the range \eqref{eq:p-pathological}.
The picture shows that the algorithm is capable of separating the two regimes correctly. 
The blue line corresponds to the line $p^* = p_{\rm true}$, and the distance on the $x$-axis between
each point and the blue line shows the absolute reconstruction error.
\review{Both results for noise-free and noisy data (10\% of intensity) are shown).
In particular,  the largest error obtained among the normal-pressure scenario is of 300 dyn/cm\textsuperscript{2} 
(less than 2\% of of the ventricular pressure in that particular case), while the largest error among the simulations with
increased pressures is of 370 dyn/cm\textsuperscript{2} (about 3\%).}
}
\label{fig:p_ven_rec}
\end{figure}

The results for the ventricular pressure reconstruction simulations are shown in Figure \ref{fig:p_ven_rec}.
As such, ventricular pressure can be evidently classified between normal and increased levels of pressure based on the reconstructed values.
Furthermore, the illustration confirms the accuracy of the reconstruction algorithm by comparing the estimated ventricular pressure with the one computed from the true ventricular pressure solution, $p_{\rm true}$.
\rreview{The relative errors between ground truth and reconstructed pressure difference
are below 2\% for the scenarios with normal pressure levels, and about 3\% for the scenarios with increased pressure levels, respectively.
Moreover, the results do not differ sensibly adding Gaussian noise to the data, since the presence of noise
does not alter the qualitative classification of pressure differences as \textit{physiological} or \textit{increased}.}

\section{Conclusions}\label{sec:conclusions}

This contribution presents a novel in silico procedure of the Parametrized-Background Data-Weak (PBDW) method for the assimilation of medical imaging data.
In this contribution, we focused on assimilating pertinent displacement data produced through magnetic resonance elastography of the human brain.
The main contribution of the proposed approach is the possibility of characterizing non-invasively intracranial pressure from partially available displacement data.
PBDW data assimilation is used to reconstruct, from the available observations, displacement and pressure fields in the whole domain.
The procedure is based on anatomical images -- used to generate a personalized computational model of the brain organ -- a finite element method for solving the underlying tissue poroelastic mechanics and a range for the values of the mechanical parameters.
From the reconstructed pressure field, it is then possible to extract clinically relevant biomarkers.
The method naturally handles partially available data and uncertainty in the physical parameters, and it can therefore overcome challenges due to missing information (e.g., pressure boundary conditions) and cope with image acquisition constraints (e.g., location of data, typically limited to one or few slices in the case of elastography). 
Moreover, since the physical solution is reconstructed in the whole domain, the proposed approach can also estimate quantity of interests when the region of interest is not observed directly.

The accuracy of the data assimilation for the joint state reconstruction (for both displacement and pressure) has been investigated in detail.
The algorithm has been validated considering 18 cases sampled within a parameter space that includes mechanical parameters (Young modulus, Poisson modulus, tissue permeability) and CSF pressure.
The simulation experiments demonstrated that including the CSF ventricular pressure in the training step, in which the manifold and the corresponding reduced space are computed, is capable of identifying patients with increased ventricular pressure based solely on partially available displacement data.

To validate the data assimilation framework, this work considers synthetic displacement data. 
However, the elastography data acquisition is part of a recent protocol used in medical imaging research in combination 
with inversion recovery (IR) protocols for identifying
tissue porosity maps \cite{lilaj_etal_2021a, lilaj_etal_2021b},
\review{differentiating between white and grey matter. This step will allow obtaining  preliminary estimates 
of different brain regions, considering the variability of the relevant mechanical and physical parameters.
}
The validation for the available in vivo data is planned for follow-up work.
However, further improvements of the proposed data assimilation algorithm that will be addressed in future research efforts include using nonlinear tissue biomechanics (e.g., viscoelasticity), anisotropic constitutive models, or enhancing the modeling through multiscale formulations (e.g., \cite{heltai-caiazzo-2019, heltai_caiazzo_mueller_2021}).
Also, other challenging extensions of the presented research concern the possibility of tackling the reconstruction of the ICP gradient instead of the whole pressure field.
Including measurements from tomo-elastography \cite{Shahryari-19} or ultrasound elastography \cite{kreft_etal_2022} in the training phase will considerably improve the proposed data assimilation approach and help towards clinical translation.

\subsection*{Acknowledgments}

This research is funded by the Deutsche Forschungsgemeinschaft (DFG, German Research Foundation) under Germany’s Excellence Strategy - {MATH+}: The Berlin Mathematics Research Center [EXC-2046/1 - project ID: 390685689].
V.V. wishes to acknowledge the financial support of the Cyprus Cancer Research Institute through the Bridges in research excellence CCRI\_2020\-\_FUN\_001, Project ``PROTEAS'' [grant ID: CCRI\_2021\_FA\_LE\_105].

\bibliographystyle{unsrt}
\bibliography{literature}

\end{document}